\RequirePackage{fix-cm}

\documentclass[smallcondensed]{svjour3}     

\smartqed  

\usepackage{mathptmx}      
\usepackage{color}

\usepackage[left=1.0in,right=1.0in,top=1in,bottom=1in]{geometry}
\usepackage{graphicx}
\usepackage{amsmath}
\usepackage{amssymb}
\usepackage{amsfonts} 
\usepackage{amsxtra} 
\usepackage{upgreek}
\usepackage[colorlinks=true,linkcolor=blue]{hyperref} 
\usepackage{tabularx} 
\usepackage{wasysym} 
\usepackage{mathtools}
\usepackage{enumerate}
\usepackage{textcomp}
\usepackage{mathrsfs}
\usepackage[section]{placeins}
\usepackage[most]{tcolorbox}

\newtcolorbox{mytextbox}[1][]{%
  sharp corners,
  enhanced,
  colback=white,
  height=10cm,
  attach title to upper,
  #1
}
\usepackage{multirow} 
\usepackage{isomath} 
\usepackage{txfonts} 
\usepackage{tensor} 
\usepackage{pifont}
\usepackage{bm} 
\usepackage{framed}
\usepackage{placeins}

\usepackage{ulem}

\usepackage{caption}
\usepackage{subcaption}

\newenvironment{frcseries}{\fontfamily{frc}\selectfont}{}

\title{Deep Learning Discrete Calculus (DLDC): A
Family of Discrete Numerical Methods by
Universal Approximation for STEM Education to Frontier Research}

\author{Sourav Saha\textsuperscript {1} \and
        Chanwook Park\textsuperscript {2} \and
        Stefan Knapik\textsuperscript {2} \and
        Jiachen Guo\textsuperscript {1} \and
        Owen Huang\textsuperscript {3} \and
        Wing Kam Liu\textsuperscript {1,2}{\dag}
}
\institute{
    {\textsuperscript{1} Theoretical and Applied Mechanics, Northwestern University, 2145 Sheridan Rd, Evanston, 60208, Illinois, USA.\\
    \textsuperscript{2} Department of Mechanical Engineering, Northwestern University, 2145 Sheridan Rd, Evanston, 60208, Illinois, USA.\\
    \textsuperscript{3} Department of Mathematics, Princeton University, Washington Road, Princeton, 08540, New Jersey, USA.\\
    \textsuperscript{\dag}Corresponding Author: Wing Kam Liu, 
	\email{w-liu@northwestern.edu}           
}
}

\begin{document}
	
\maketitle
\begin{abstract}
The article proposes formulating and codifying a set of applied numerical methods, coined as \textbf{D}eep \textbf{L}earning \textbf{D}iscrete \textbf{C}alculus (DLDC), that uses the knowledge from discrete numerical methods to interpret the deep learning algorithms through the lens of applied mathematics. The DLDC methods aim to leverage the flexibility and ever increasing resources of deep learning and rich literature of numerical analysis to formulate a general class of numerical method that can directly use data with uncertainty to predict the behavior of an unknown system as well as elevate the speed and accuracy of numerical solution of the governing equations for known systems. The article is structured in two major sections. In the first section, the building blocks of the DLDC methods are presented and deep learning structures analogous to traditional numerical methods such as finite difference and finite element methods are constructed with a view to incorporate these techniques in Science, Technology, Engineering, Mathematics (STEM) syllabus for K-12 students. The second section builds upon the building blocks of the previous discussion, and proposes new solution schemes for differential and integral equations
pertinent to multiscale mechanics. Each section is accompanied with mathematical formulation of the numerical methods, analogous DLDC formulation, and suitable examples.  

\end{abstract}


\section{Introduction}
The problems of engineering and physical science can be categorized into three types \cite{saha2021hierarchical}, type 1 problems are the problems with limited physical understanding and a lot of experimental data, type 2 problems are those with incomplete physical understanding and some experimental data, and type 3 problems are those for which there is sufficient knowledge about the system but solving for the system response is computationally challenging. An example of type 1 problem is relating the spatio-temporal variation of temperature to the resulting mechanical properties in metal additive manufacturing \cite{xie2021mechanistic}, a type 2 problem would be calibrating the heat source and other relevant models for computational fluid dynamics simulation for metal additive manufacturing process \cite{gan2021benchmark}, and a type 3 problem would be relating the microstructure and fatigue life of a 3D printed metallic part \cite{xie2021mechanistic,kafka2021image}. Although the classification is not a strict one, it points to two facts: a) even with numerical modeling, aid of data is needed at some level, and b) by combining data science and numerical method one can cover the spectrum of solving for completely unknown phenomenon to challenging problems. Hence, the question naturally arises, is there a way we can come up with a new method that gives us the best of both worlds? 

Establishing this connection between numerical methods and deep learning is the most sought after prize in the literature nowadays. The quest started with solving the dynamic system directly using data science and deep learning methods \cite{sirignano2018dgm,raissi2018deep,han2017deep}. The deep learning methods primarily include deep neural network \cite{yu2018deep}, recurrent neural network \cite{xiao2018nonlinear}, convolutional neural network \cite{zhu2018convolutional}, and residual neural network \cite{kani2017dr}. Other techniques from data science such as unsupervised learning are also heavily used to analyze the data and extract meaningful features that may or may not have explicit physical meaning (sometimes called the latent variables) \cite{karl2016deep,chakraverty2017artificial}. However, it was apparent that using directly deep learning methods results in lack of generalization \cite{magill2018neural}. The neural networks are essentially highly non-linear interpolation functions with parameters being trained on observations \cite{liu2021mechanistic}. Hence, to predict the behaviour of a system outside the training range becomes challenging. Moreover, the neural networks in such methods may become very complex and it is often impossible to interpret the inner working of the neural network in use. This results in a "black box" method which we can just use without proper know-how. Finally, these methods are data-hungry \textit{i.e.,} it requires a lot of data to implement these.  

As a solution to the lack of generalization, combining the knowledge of mechanics with data science is proposed as \textit{Mechanistic Data Science} (MDS) \cite{liu2021mechanistic}. The knowledge of mechanics can be incorporated in selection of features and curing of data to be used as features, or, governing principles such as conservation of energy can be directly used as an optimization constraint. The later method is also termed as the Physics Informed Neural Network (PINN) \cite{raissi2019physics,karniadakis2021physics,pang2019fpinns}. These methods have become extremely popular as these serve perfectly to solve type 3 problems \cite{krishnapriyan2021characterizing,li2021physics,leung2021nh}. Moreover, these methods only partially solve the issue of the requirement of large datasets. However, neither of these theories solves the problem of lack of interpretation as the neural network is still there. Finally, for type 1 problems (where we only have data) the idea of adding governing principle does not work.

Data science techniques to discover the underlying governing equation from data have been proposed \cite{brunton2016discovering,kaiser2018sparse,kaheman2020sindy}. Another example is Dimension-Net which discovers the underlying non-dimensional numbers directly from data \cite{gan2021universal,xie2021data}. While these techniques are exciting, one needs to have some idea about the system behavior to use these techniques as these are primarily regression techniques working on a library of candidate functions. Therefore, there is still a need for a data science framework that can directly use the experimental data to solve for the dynamics of the system and at the same time interpret the technique. 

Efforts towards merging numerical methods with deep learning have gained increased attention in the last few years. The initial efforts were to establish the parallels between a subset of numerical methods with a subset of deep learning techniques \cite{chen2018neural,zhang2021hierarchical}. For example, how a recurrent neural network behaves like an ordinary differential equation \cite{niu2019recurrent} such as wave equation \cite{hughes2019wave}, or, how some finite difference techniques or multi-grid methods have parallels with convolutional neural network \cite{he2019mgnet}. While these works are important as an interpretation of the deep learning techniques, they do not suggest how to improve the numerical solvers using data science or vice-versa. One of many early efforts that tries to leverage deep learning to improve existing numerical methods is called Hierarchical Deep Learning Finite Element Method (HiDeNN-FEM) \cite{zhang2021hierarchical} where using custom neural networks different linear and non-linear interpolation shape functions are obtained. A recent breakthrough in interpretable general deep learning methods for solving partial differential equation came in the form of operator learning \cite{kovachki2021neural}. Two major methods were proposed concurrently for operator learning called Fourier Neural Operator (FNO) \cite{li2020fourier}, and DeepONet \cite{lu2019deeponet}. There have been studies showing comparisons between these two techniques \cite{lu2022comprehensive}. Keeping the similarities and differences aside, the philosophy of both of these methods are same: mapping functions from input space to output space instead of mapping data. Based on the same philosophy, several other structures are being proposed including graph kernel network and non-local kernel network \cite{anandkumar2020neural,you2022nonlocal}. 
\begin{figure}
    \centering
    \includegraphics[width=\linewidth]{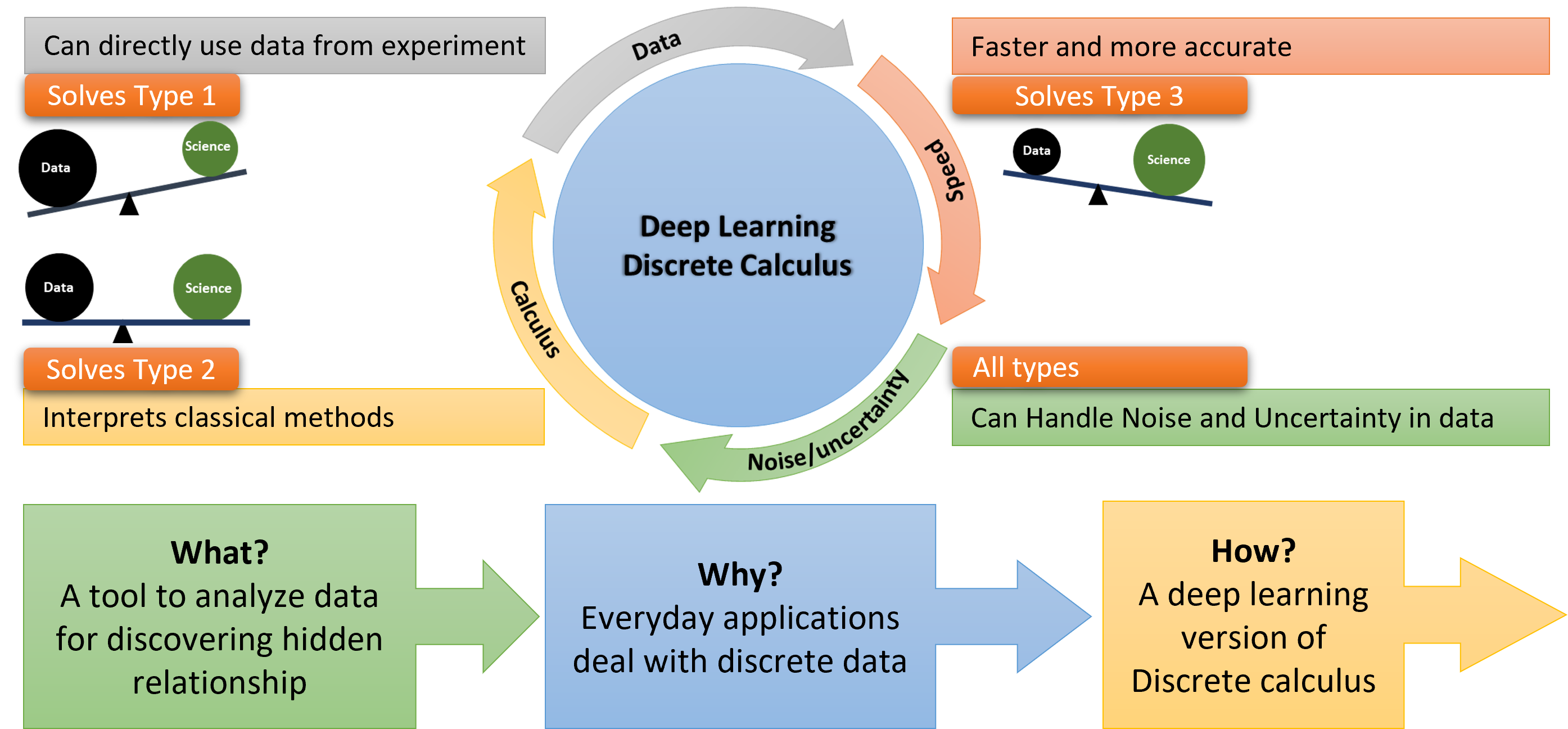}
    \caption{The key features of the Deep Learning Discrete Calculus (DLDC) method.}
    \label{fig1}
\end{figure}

Aside from solving challenging problems in scientific and engineering research, contemporary deep learning methods have a huge pedagogical potential as well. In order to develop a competent workforce, big tech companies and entrepreneurs are now focusing on integrating artificial intelligence (AI) in high school curriculum \cite{tucker2020exploring}. There is increasing need for skilled workforce over outstanding academic achievers in the industries. Naturally, the question arises how to train the young students in AI and at the same time, make the curriculum interesting for them.  If the AI is introduced in a traditional way through teaching computer science, the diversity in the student body will become very difficult to attain. Due to the socioeconomic background of the students, direct training in machine learning algorithms may become very hard to attain. Thus untested introduction to AI for middle and high school students will create another bias. On the other end of the spectrum, the way mathematics, especially, calculus is taught in the high schools may become too abstract for some students to grasp. It is natural for the young students to be interested to learn a concept with a clear idea about how to implement it. Usually, in the school curriculum the calculus is taught through some textbook examples which often make the students disinterested in the topic \cite{wang2022perspective}. This lack of interest results in reduced number of students choosing Science, Technology, Engineering, and Mathematics (STEM) majors. There have been a lot research \cite{touretzky2022artificial,touretzky2022machine,yin2022ai4all} on how to incorporate data science and AI into high school education. However, most of these studies are at strategic level and do not provide a clear pathway to merge two broad fields of calculus and AI. It is the authors' opinion that introducing AI through deep learning and interpretation of deep learning via elementary calculus is the most optimized way to teach calculus at the high school level. This has two-way benefit: one, implementing calculus through advanced computational tools with real-life data will be possible for young students. It will increase their interest and understanding on the basic calculus. Two, interpretability of deep learning algorithms can be improved if the algorithms are thought from the first principle.        

With the broad scope in mind, in this article, the authors propose an applied numerical method that aims to combine the mechanistic knowledge, classical numerical methods, and deep learning algorithms called Deep Learning Discrete Calculus (DLDC) (see, Figure \ref{fig1}). A formal definition of Deep Learning Discrete Calculus can be given as:

\begin{tcolorbox}[colback=white!5!white,colframe=black!75!black]
 The Deep Learning Discrete Calculus is a discipline to integrate fundamental definitions of calculus with numerical methods using deep learning neural network.

\end{tcolorbox}

Deep learning is an extremely powerful tool for feature extraction and prediction, but it requires massive and cleansed data sets to perform well. For many engineering applications, the collection of so much data is prohibitively expensive or even impossible, making traditional deep learning approaches nonviable. Incorporating discrete calculus and numerical method tools into neural network architectures can empower deep learning to address engineering problems where datasets may be small and noisy. Additionally, numerical method-informed deep learning frameworks sometimes have the potential to extrapolate beyond available observations, which is something traditional deep learning implementations are infamously incapable of. The synthesis of numerical analysis ideas with machine learning can reduce the expenses of collecting data, pre-processing data, and training models, while improving the interpretability, accuracy, and generalizability of deep neural networks. Figure \ref{fig1} shows the key features of the DLDC method. The DLDC aims to become a general method for solving the three different types of problems discussed before. 
The article has two major parts. The first part aims to introduce the concept of DLDC through building blocks of calculus, such as, differential and integral calculus for STEM education. In later sections, the article discusses how to apply the DLDC methods to solve differential and integral equations, and how it has advantages over the traditional numerical methods. Each demonstration is accompanied by an example to further reinforce the ideas presented in the article.

\section{Formulating DLDC for STEM Education}\label{sec2}

\begin{figure}[h]
\centering
\includegraphics[width=0.85\textwidth]{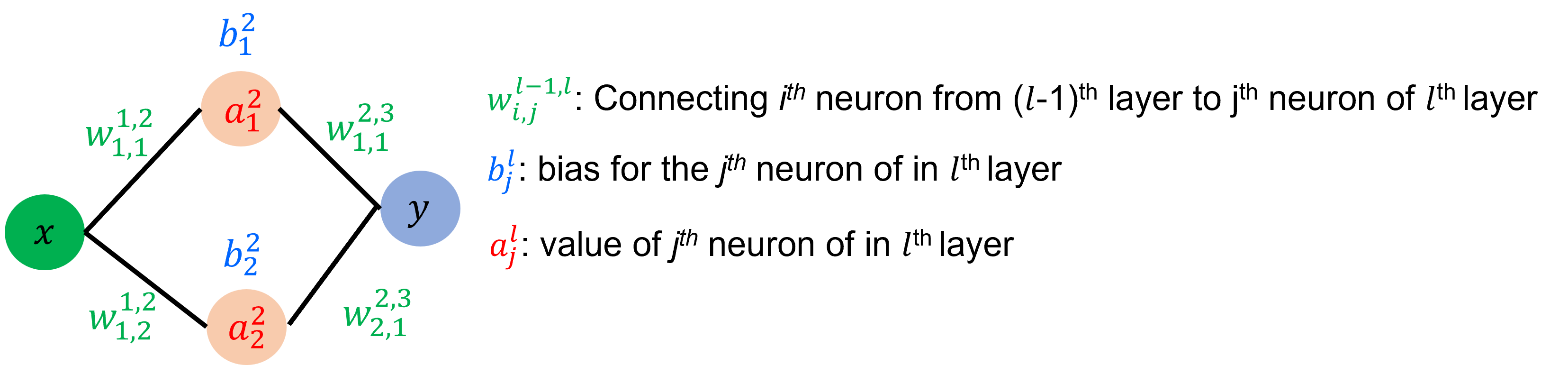}
\caption{A neural network model with input, hidden, and output layers.}
\label{fig:FigNN}
\end{figure}

In this section, the basic mathematical foundations of the DLDC methods are presented. This section will be particularly useful for implementing in STEM education. Before going into too much detail, a brief discussion on neural network is required. For extensive discussion, the readers can refer to \cite{liu2021mechanistic}. A neural network has three types of layers, input, hidden, and output. The input variables of interest go into input layer, and the output variables are obtained through the output layer. The hidden layer(s) takes input variable and pass those through a non-linear function, called activation function, with optimization parameters. Mathematically, if $x$ is the input variable, $\mathcal{A}$ is the activation function, $y$ is the output variable, the equation for a neural network can be written as,

\begin{equation}
\centering
\label{Eq1}
y =  \mathcal{A}(Wx+b)
\end{equation}
where $W$ is the weight, and $b$ is the bias. By varying the weight and bias to minimize the error between predicted output and training data, one can achieve an approximation of any linear and non-linear functional relationship. A schematic relationship showing the details of a neural network with one input, one hidden, and one output layer is shown in Figure \ref{fig:FigNN}. The convention of super and sub-scripts for weights and biases are shown in Figure 2. Following this convention, the output can be written as a function of the trainable weights, biases, and given input as

\begin{equation}
\centering
\label{Eqpp}
y =  \mathcal{A}\left[W_{1,1}^{2,3}(\mathcal{A}(W_{1,1}^{1,2}x+b_{1}^{2}))+W_{2,1}^{2,3}\mathcal{A}(W_{1,2}^{1,2}x+b_{2}^{2})\right]
\end{equation}

The trainable $\bm{W}$ and $\bm{b}$ parameters are varied to minimize the following loss function:

\begin{equation}
    Loss function = \frac{1}{P}\sum_{n=1}^{P}(y-y^{*})^2
\end{equation}

\subsection{Differential Calculus}
Based on this basic definition of neural network, we proceed to show how to build a neural network from scratch to mimic the finite difference methods. Let us consider a function $f(x)$ as shown in Figure 2(a). If we want to have an approximation of first derivative of $f(x)$ at A, the forward difference method will give us \cite{chapra2011numerical}, 
\begin{equation}
\centering
\label{Eq2}
\frac{dy}{dx} =  \frac{f(x_{j+1})-f(x_j)}{x_{j+1}-x_{j}}
\end{equation}
Here, the index $j$ indicates point A, and $j+1$ indicates the next point on the function. The central difference method gives, 

\begin{equation}
\centering
\label{Eq3}
\frac{dy}{dx} =  \frac{f(x_{j+1})-f(x_{j-1})}{x_{j+1}-x_{j-1}}
\end{equation}
A higher order approximation of the first derivative involving three points looks like,

\begin{equation}
\centering
\label{Eq4}
\frac{dy}{dx} =  \frac{-3f(x_j)+4f(x_{j+1})-f(x_{j+2})}{x_{j+1}-x_{j-1}}
\end{equation}
Therefore, from observation, we can simply write the general form of these equations as, 

\begin{equation}
\centering
\label{Eq5}
\frac{dy}{dx} =  \sum_{j=1}^{n}w_{j}f(x_j)
\end{equation}

\begin{figure}[h]
\centering
\includegraphics[width=\textwidth]{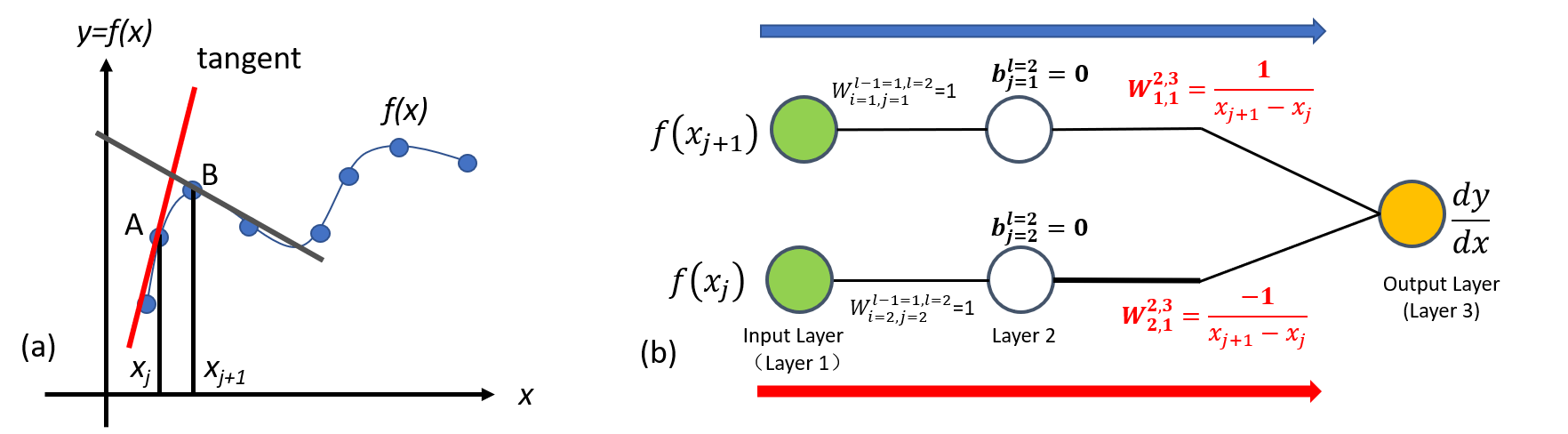}
\caption{(a) A schematic diagram showing a continuous function $f(x)$ and discrete sampling data, (b) The DLDC structure for the first derivative (forward difference method).}
\label{fig:Fig2}
\end{figure}

This form for numerical differentiation is famously known as differential quadrature method in numerical analysis \cite{bellman1971differential,bert1996differential,shu2012differential}. If we compare the form of Eqn. \ref{Eq5} with Eqn. \ref{Eq1}, a clear parallel can be observed. The finite difference formulas can be obtained using linear activation functions with zero bias. Inspired from this observation, a first-order derivative equivalent neural network is proposed in Figure 2(b). In the figure, the weights and biases are expressed in index notation for convenience. In $W_{i,j}^{l-1,l}$, the subscripts $i$ and $j$ denotes the incoming neuron and outgoing neuron, respectively. The superscript $l$ identifies the layers. The index $l-1$ indicates the incoming layer and $l$ indicates the outgoing layer. The network has one hidden layer, rectified linear unit activation function, and zero bias. The weights between input and hidden layer are constrained to be 1, and the weights between the hidden and output layer (red marked) are optimized. All the activation functions are linear and biases are set to zero. The cost function can be written as, 

\begin{equation}
\centering
\label{Eq6}
\frac{1}{N} \sum_{i=1}^{N}\left(\frac{dy}{dx}^{true}-\frac{dy}{dx}^{prediction}\right)^{2}_{i}
\end{equation}

Here, $N$ is the sample size. It is interesting to note that after optimization, the trainable weights becomes $W_{1,1}^{2,3}=\frac{1}{x_{j+1}-x{j}}$ and $W_{2,3}^{2,1}=\frac{-1}{x_{j+1}-x{j}}$. This value resembles the forward difference method. If a neural network like this is trained once against data generated by a known function and it's derivative, the trained network can be used to predict the first derivative of any function provided the sampling remains the same. 

\begin{figure}[h]
\centering
\includegraphics[width=0.95\textwidth]{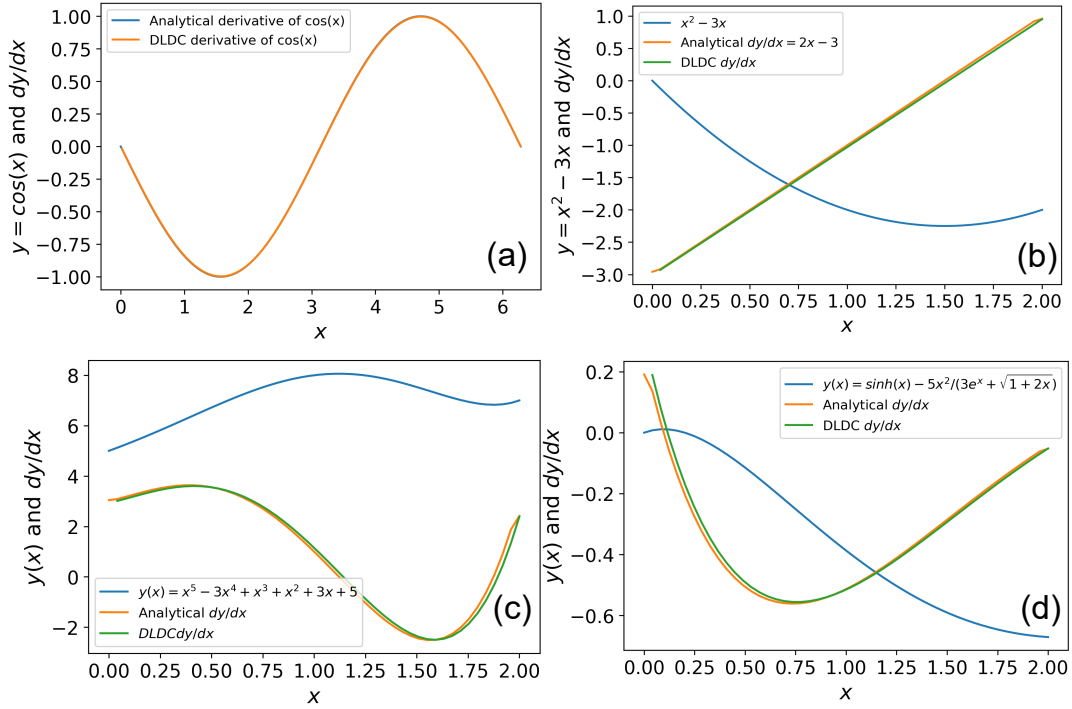}
\caption{The DLDC results for first-order derivative. (a) Performance for predicting the first derivative of cosine function. The training was done with sine function. (b) Training performance for simple polynomial, and (c), (d) performance of the network to predict derivatives of more complicated functions.}
\label{fig:DLDC_nodes}
\end{figure}

This concept has been applied to formulate the first order derivative of different types of known function. Some of the results are presented in Figure \ref{fig:DLDC_nodes}. In case of Figure \ref{fig:DLDC_nodes}(a), the neural network is trained for a simple sinusoidal function in a range between $0-360$ degrees. With the trained neural network, we tried to predict the first-order derivative of cosine function which gives very good prediction. In the next example, see, Figure \ref{fig:DLDC_nodes}(b), a simple polynomial was used to generate training data for the neural network. Once the network is trained, the same network is applied to predict the derivative of a higher-order polynomial (Figure \ref{fig:DLDC_nodes}(c)) and a relatively complicated function (Figure \ref{fig:DLDC_nodes}(d)). 

\subsection{Integral Calculus}
\begin{figure}[h]
\centering
\includegraphics[width=1\textwidth]{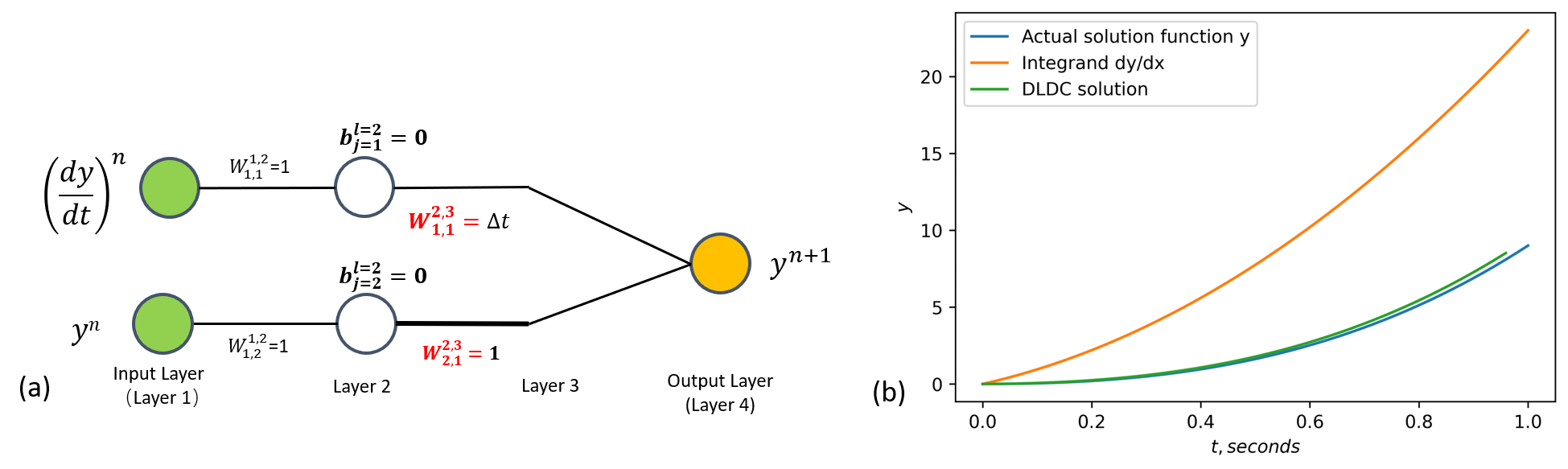}
\caption{(a) The DLDC structure for the explicit Euler formulation, (b) A comparison between the analytical and the DLDC prediction for the time integration.}
\label{fig:Fig4}
\end{figure}

The next building block is to solve an integration problem. Numerically, an integral of a function with respect to a variable is an estimation of area under the curve representing that function when plotted against the said variable. Based on this simple concept, numerical methods such as trapezoidal or Simpson's method are in use. However, engineers are mostly interested to solve time dependent ordinary or partial differential equation. The first method that comes to mind for such computation is the Euler's method. Let us suppose we want to solve an ordinary differential equation.    

\begin{equation}
\centering
\label{Eq7}
\frac{dy}{dt} = f(y,t)
\end{equation}
where $y$ is the dependent variable and $t$ is time. A DLDC equivalent structure for explicit Euler method and associated prediction is shown in Figure \ref{fig:Fig4}. The generalized alpha method for Euler's method is: 

\begin{equation}
\centering
\label{Eq8}
y^{n+1} = y^n + \Delta t \left([1-\alpha].f^n+[\alpha].f^{n+1}\right),
\end{equation}
where, $n$ is the time step, $\Delta t$ is step size, and $\alpha$ is a factor that indicates an explicit method if the value is 0 and implicit method if the value is 1. The structure of the DLDC for Euler integration is exactly the same as forward difference method, except the inputs are the values of the integrand and dependent variable at the previous step and the output is the value of the dependent variable in the next time step. In the neural network shown in Fig. \ref{fig:Fig4} (a), the biases are set to zero and activation functions are linear. The black weights are fixed to a value of 1 and the red weights are only allowed to train. The values of the red marked weights in the figure are the final values of these weights after training. It turns out that after the training, the weights finally get to a value that makes the output mimic the equation Eq. \ref{Eq8}.

The numerical prediction in the example was made by solving $\frac{dy}{dt}=15t^2+8t$ with $y=0$ initial condition. The predicted and analytical solutions diverge slightly near 1 second. But this error is coming from the Euler method itself as the weights of the trained DLDC is exactly same as the explicit Euler method. Similar network can be constructed for implicit Euler method as well.

\begin{figure}[h]
\centering
\includegraphics[width=0.65\textwidth]{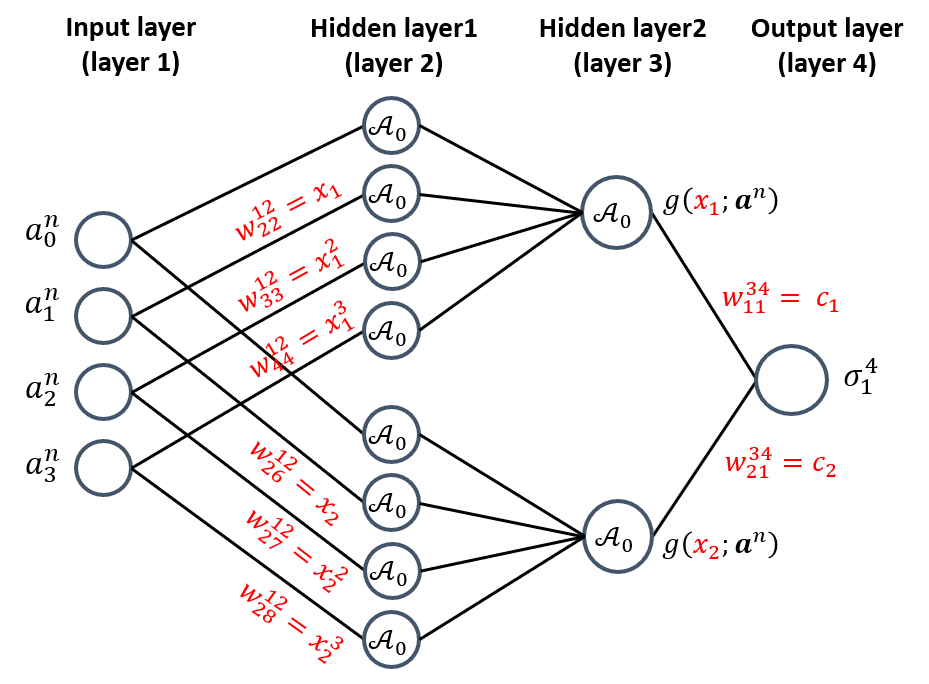}
\caption{The DLDC structure for 2-point Gauss quadrature method. The red marked terms are to be optimized during training.}
\label{fig:Fig5}
\end{figure}

\begin{figure}[h]
\centering
\includegraphics[width=0.85\textwidth]{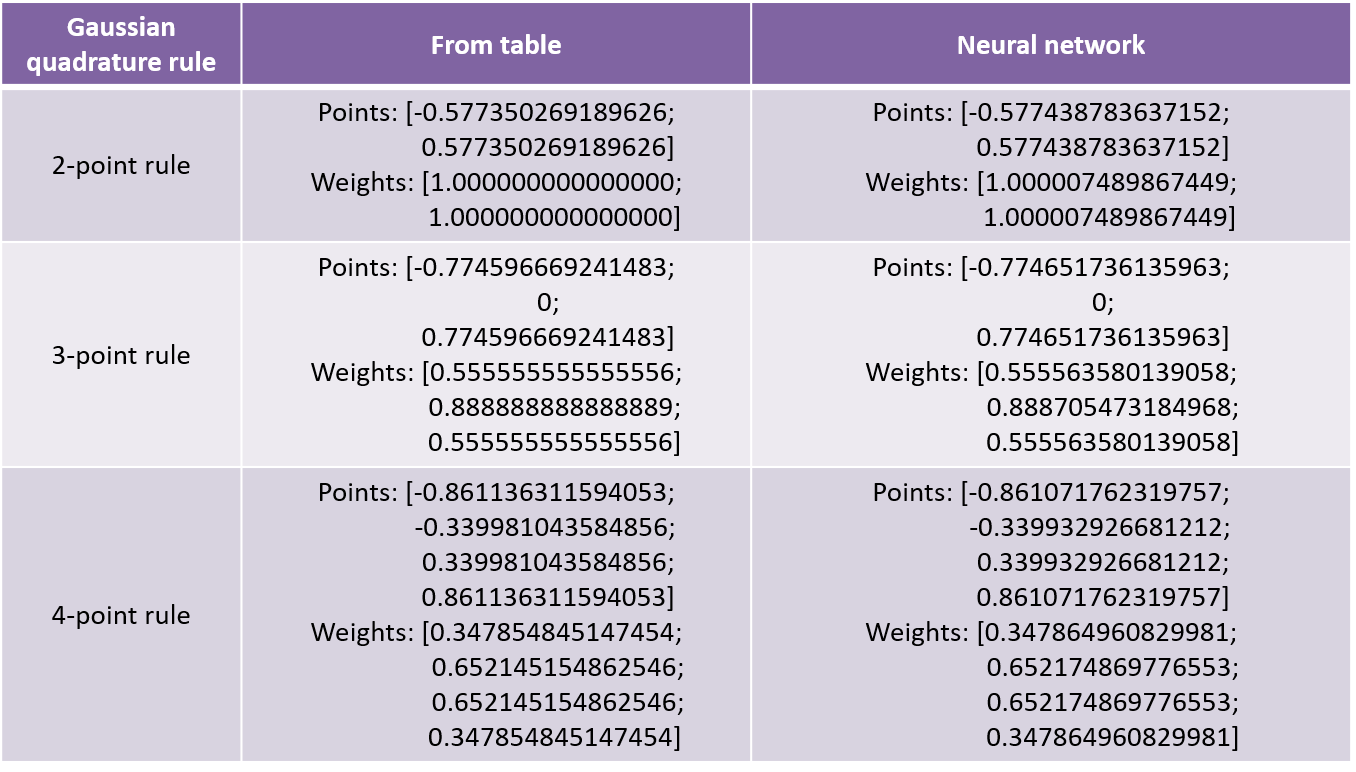}
\caption{A table showing the comparison between different quadrature methods and optimized values for the DLDC.}
\label{fig:Fig6}
\end{figure}

The next attempt was made to construct a DLDC method for Gauss quadrature integration method \cite{chapra2011numerical}. In brief, the Gauss quadrature method can perform integration by converting integral to fixed limit between -1 to 1. The equations look like, 

\begin{equation}
\centering
\label{Eq9}
\int_{a}^{b} f(x)dx = \int_{-1}^{1} f\left(\frac{b-a}{2}\xi+\frac{a+b}{2}\right)\frac{dx}{d\xi}d\xi,
\end{equation}

\begin{equation}
\centering
\label{Eq10}
\int_{-1}^{1}g(x)dx = \sum_{i=1}^{n}w_{i}g(x_{i}).
\end{equation}
Here, $n$ is the number of quadrature points and $w_{i}$ are weights assigned to each point. $g(x_{i})$ is the value of the converted function at each quadrature points. The DLDC structure for the 2-point Gauss quadrature is shown in Figure \ref{fig:Fig5}. In order to train the network, a third order polynomial $g(x,\textbf{a})=a_0+a_1 x+a_2 x^2+a_3 x^3$ is assumed and the corresponding Gauss integral $\int_{-1}^{1}g(x,\textbf{a})dx=2a_0 + \frac{2}{3}a_2$. The training data is generated for different values of $\textbf{a}$ and the integral. In the DLDC neural network, the activation function is rectified linear unit between input and first hidden layer and at the last layer, and linear for the second hidden layer. If we train the neural network with respect to $x_1 , x_2, w_{11}^{34}, w_{21}^{34}$ with following cost function, 

\begin{equation}
\centering
\label{Eq11}
\frac{1}{N}\sum_{n=1}^{N}\left(\sigma_{1}^{4}(x_1 , x_2 , w_{11}^{34}, w_{21}^{34};\textbf{a}^n)-\int_{-1}^{1}g(x,\textbf{a})dx\right)^2,
\end{equation}

where 
\begin{equation}
\centering
\label{Eq12}
\sigma_{1}^{4}(x_1 , x_2 , w_{11}^{34}, w_{21}^{34};\textbf{a}^n)=\sum_{j=1}^{N_{i-1}}w_{jk}^{i-1,i}\sigma_{j}^{i-1}+b_k^{i}=c_1 g(x_1 ; \textbf{a}^n) + c_2 g(x_2, \textbf{a}^n),
\end{equation}
after the training the quadrature coordinates $x_1$, $x_2$, and the weights $c_1$ and $c_2$ have exact same value as theoretical Gauss quadrature method. The structure was extended to 3 and 4-point rules and for all the cases the trained values for the co-ordinates and weights match well with the analytical value (see Figure \ref{fig:Fig6}). Just like the Euler method and finite difference approximation, the Gauss quadrature can be trained with a known integration can be used for any other integration without retraining. 

\section{DLDC Formulation for Frontier Research} 

This section builds upon the ideas of the DLDC formulation for STEM education and proposes discrete calculus-inspired deep learning methods to solve partial differential equations and integral equations which are at the core of applied mechanics. 

\subsection{Solving Partial Differential Equations}

\subsubsection{Convolution Hierarchical Deep-learning Neural Network (C-HiDeNN)}

In this section, we use DLDC to formulate a new interpolation theory based on convolution, called Convolution Hierarchical Deep-learning Neural Network (C-HiDeNN). C-HiDeNN interpolants are smoother and more accurate than those of FEM. They are designed to have arbitrary reproducing polynomial orders, resulting in  superior convergence behavior without using higher order elements.

To understand C-HiDeNN, we should review HiDeNN-FEM first \cite{saha2021hierarchical,zhang2021hierarchical}. HiDeNN-FEM is a neural network representation of FEM. HiDeNN-FEM shape functions are small neural networks that are designed to be mathematically equivalent to FEM shape functions. Therefore, the weight and bias of HiDeNN-FEM shape functions are represented as a function of nodal coordinates and nodal field variables. The element-wise shape functions are then gathered to form a global neural network and finally, a loss function based on the principle of minimum potential energy is minimized. 

In HiDeNN-, the solution is obtained by optimizing nodal coordinates and nodal field variables with respect to the potential energy loss function. There are two solution schemes for the HiDeNN-FEM. First, the nodal coordinates are kept fixed and only the nodal variables are updated. This makes the HiDeNN-FEM mathematically equivalent to FEM, thus the solution accuracy and computation time of HiDeNN-FEM is on the same order as FEM. Second, both the nodal coordinates and nodal variables are updated, thus making it equivalent to r-adaptive FEM. For detailed discussions, readers may refer to \cite{zhang2021hierarchical}.

\begin{figure}[h]
\centering
\includegraphics[width=0.95\textwidth]{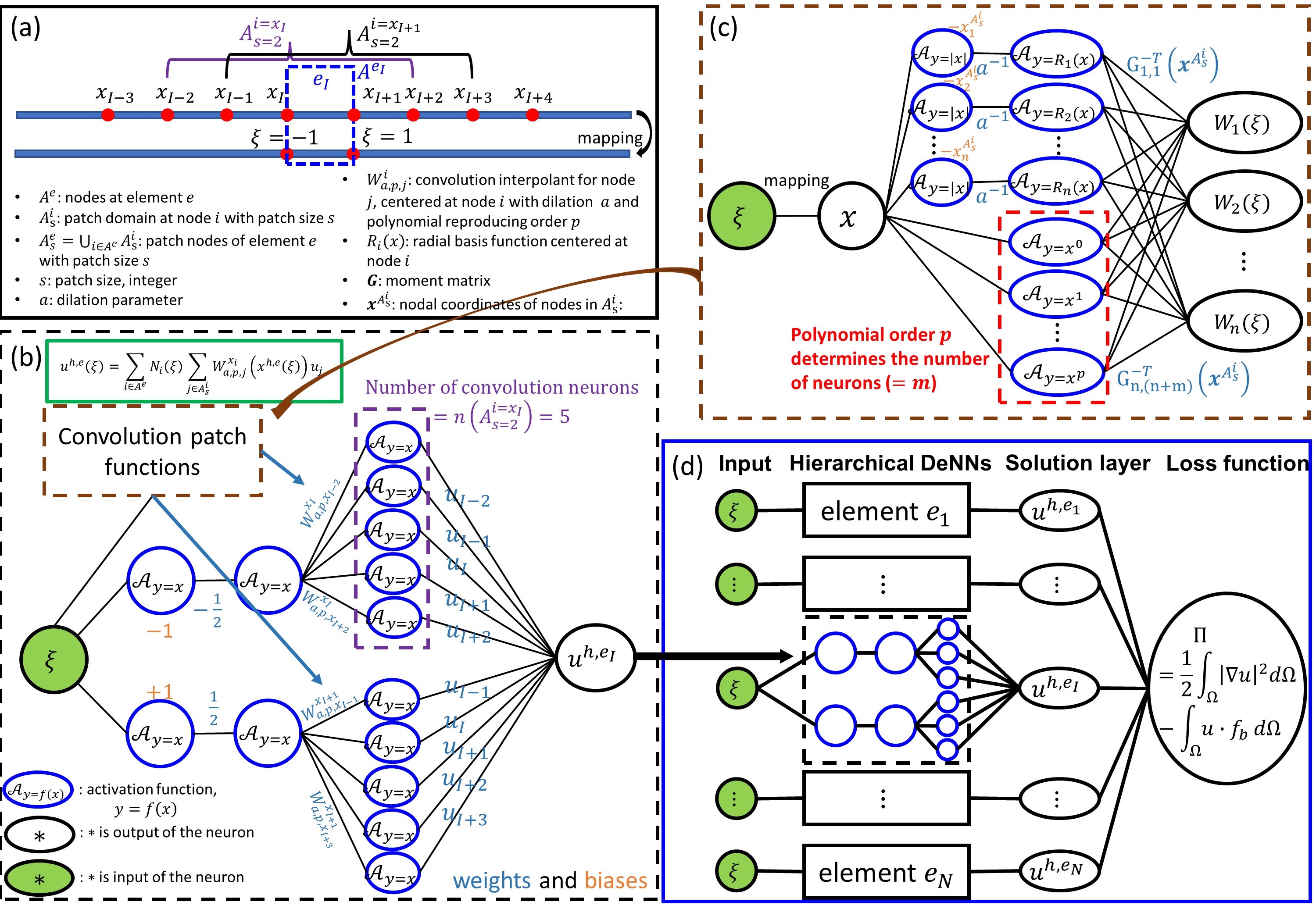}
\caption{C-HiDeNN formulation for 1D Poisson problem. Light blue and orange terms are weights and biases of the neural network, respectively. If there is no weight or bias assigned for a neuron, it will have fixed negligible weight=1 and bias=0. Functions inside neurons with blue edges represent activation functions while those with black edges represent inputs (green color) and outputs (white color) of the neuron. (a) shows nodal coordinates in both physical and parametric space, focusing on the element of interest $e_I$. Nodal patch domains and other terminologies are defined below. (b) represents the C-HiDeNN shape function of element $e_I$, which constitutes the hierarchical DeNN layer of the global neural network (d). (c) is the convolution patch function that can be found at the greed dotted box in (b). This figure is taken from \cite{park2023convolution}.}
\label{fig:C-HiDeNN formulation}
\end{figure}

Since HiDeNN-FEM is still based on finite element interpolations, the improvement of global solution accuracy is insignificant even with the automatic r-adaptivity. This prompted the idea of convolution to answer the question: "Can we incorporate convolution filters to HiDeNN-FEM shape functions to achieve highly smooth and accurate interpolants?". This idea is written as (in 1D):

\begin{equation} \label{Eq C-HiDeNN 1}
    u^{h,e}(\xi) = \sum_{i \in A^e} N_i(\xi) \sum_{j \in A^i_s} W^{x_i}_{a,p,j}\left(x^{h,e}(\xi)\right)u_j = \sum_{k \in A^e_s} \tilde{N}_k(\xi)u_k,
\end{equation}

\noindent with elementwise mapping between natural coordianates and physical coordinates:
\begin{equation} \label{Eq C-HiDeNN 2}
    x^{h,e}(\xi) = \sum_{i \in A^e} N_i(\xi)x_i.
\end{equation}

In Eq. \ref{Eq C-HiDeNN 1}, two interpolants appear: 1) general polynomial interpolants, $N_i(\xi)$, defined over element $A^e$; 2) convolution patch functions, $W^{x_i}_{a,p,j}\left(x^{h,e}(\xi)\right)$, defined over nodal patch domain $A^i_s$. The elemental patch domain is visualized in Fig.\ref{fig:C-HiDeNN formulation}(a) where the integer-valued patch size $s$ refers to the number of element layers surrounding node $i$. That is, $A^{i=x_I}_{s=2}$ contains node $x_{I-2}$ through node $x_{I+2}$ while $A^{i=x_{I+1}}_{s=2}$ contains node $x_{I-1}$ through node $x_{I+3}$.

The general polynomial interpolants $N_i(\xi)$ can be thought as any finite element shape functions that are compactly supported and satisfy Kronecker delta and partition of unity. The convolution patch functions $W^{x_i}_{a,p,j}\left(x^{h,e}(\xi)\right)$ are also compactly supported supported interpolation functions that satisfy Kronecker delta and reproducing conditions. Here, the dilation parameter $a$ determines the size of support domain and $p$ is the reproducing polynomial order. For the convolution patch functions, we borrow from well-developed meshfree techniques. In this study, the radial point interpolation method is adopted because this method returns stable interpolants that satisfy Kronecker delta property and reproducing conditions. For details, readers may refer to \cite{liu2001local,liu2005introduction}.

Finally, the double summation in Eq. \ref{Eq C-HiDeNN 1} is combined into a single summation over the elemental patch domain $A^e_s=\cup_{i \in A^e}A^i_s$, and the resulting convolution interpolants are $\tilde{N}_k(\xi)$ where $k \in A^e_s$. The convolution interpolants therefore satisfy compact supportness, partition of unity, Kronecker delta, and reproducing conditions, making it easy to apply boundary conditions.

In terms of DLDC, the C-HiDeNN shape function can be written as a neural network illustrated in Fig.\ref{fig:C-HiDeNN formulation}(b). The first two hidden layers refer to the finite element shape functions, and the third hidden layer is for the convolution patch functions. Note that the number of neurons is determined by the number of nodes in nodal patch domains. That is, larger patch size $s$ leads to larger number of neurons or higher connectivity. This is analogous to the convolutional neural network (CNN) kernels. 

In C-HiDeNN, higher reproducing polynomial order $p$ is achieved by setting large patch size $s$. Based on our preliminary study, $s \geq p$ to get stable convolution patch functions. In other words, C-HiDeNN can achieve higher order $p$ by increasing the connectivity of nodes while using the same linear elements. This is completely different in FEM where higher order elements must be used to build higher order shape functions. Therefore, the global degrees of freedom (DOFs) of C-HiDeNN are the same as linear FEM, but the bandwidth of the global stiffness matrix becomes larger due to increased nodal connectivity. 

The weights of the third hidden layer in Fig.\ref{fig:C-HiDeNN formulation}(b) are built form another sub-neural network shown in Fig.\ref{fig:C-HiDeNN formulation}(c). This follows the radial point interpolation theory and details are discussed in \cite{park2023convolution}. Finally, the elementwise shape functions are assembled to form a global neural network and potential energy loss function shown in Fig.\ref{fig:C-HiDeNN formulation}.

To demonstrate superior accuracy of C-HiDeNN compared to FEM, a 2-D Poisson problem has been solved:

\begin{equation} \label{Eq C-HiDeNN 3}
\begin{split}
    \nabla \cdot (\nabla u(\textbf{x})) + b_f(\textbf{x}) &= 0  \quad \textrm{in} \  \Omega  \\
    u &= 0  \quad \textrm{on} \  \Gamma
\end{split}
\end{equation}

\noindent where $b_f(\textbf{x})$ is the body force and a square domain $\Omega$ whose lower left corner is at $(0,0)$ and upper right cornder is at (10,10) is used. $\Gamma$ is the boundary of the domain $\Omega$. We set the analytical field variable as:

\begin{equation} \label{Eq C-HiDeNN 4}
    u(\textbf{x}) = \frac{1}{625} (x^2-10x)(y^2-10y)\left(2e^{-2((x-3)^2+(y-3)^2)} + e^{-2((x-7)^2+(y-7)^2)}\right).
\end{equation}

\noindent L2 and H1 error estimators are used:
\begin{equation} \label{Eq C-HiDeNN 5}
\begin{split}
    \lVert e \rVert_{L_2} &=  \lVert u-u^h \rVert_{L_2} = \frac{\left(\int_{\Omega}(u-u^h)^2 dx\right)^{1/2}}
    {\left(\int_{\Omega}u^2 dx\right)^{1/2}} \\
    \lVert e \rVert_{H_1} &= \lVert u-u^h \rVert_{H_1} = \frac{\left( \int_{\Omega}(u-u^h)^2 dx + \int_{\Omega} \lVert \nabla u- \nabla u^h \rVert^2_2 dx \right)^{1/2}}
    {\left( \int_{\Omega}u^2 dx + \int_{\Omega} \lVert \nabla u \rVert ^2_2 dx \right)^{1/2}} 
\end{split}
\end{equation}

\begin{figure}[h]
\centering
\includegraphics[width=0.95\textwidth]{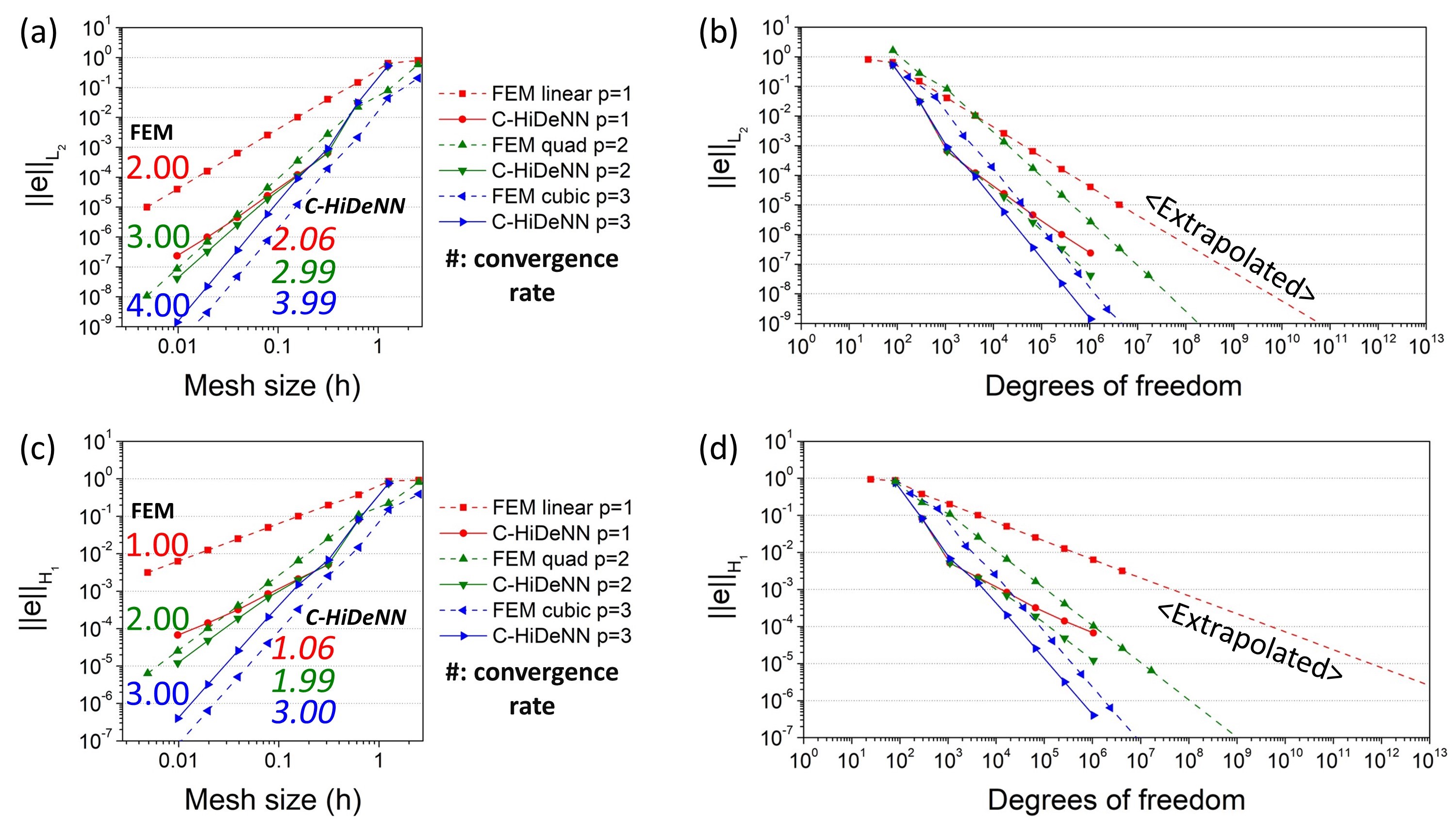}
\caption{Convergence plot. (a, b) are for $L_2$ norm error and (c, d) are for $H_1$ norm error estimation. (a, c) are the error vs. mesh size plots and (b, d) are the error vs. degrees of freedom plots. Numbers on the graph are the convergence rates (italic font for C-HiDeNN). For C-HiDeNN, patch size $s=3$ and dilation parameter $a=30$. The graphs are colored by the reproducing polynomial order $p$. FEM uses dashed lines and C-HiDeNN uses solid lines.}
\label{fig:C-HiDeNN convergence}
\end{figure}

Fig.\ref{fig:C-HiDeNN convergence} shows C-HiDeNN-FEM can have arbitrary convergence rates depending on the reproducing polynomial order $p$. That is, for the same $p$, the convergence rates (slope of Fig.\ref{fig:C-HiDeNN convergence}(a,c)) of FEM and C-HiDeNN are asymptotically the same. However, their y-intercepts are different. When $p=1$, C-HiDeNN is around two orders of magnitude more accurate than FEM for a given mesh size. The gap decreases as $p$ goes up and the FEM turns around when $p=3$ (i.e., blue dotted curves are lower than blue solid curves in Fig.\ref{fig:C-HiDeNN convergence}(a,c)). However, this is not a fair comparison as higher order FEM uses higher order elements that has more DOFs while C-HiDeNN still uses the linear elements. 

The errors are plotted with respect to DOFs in Fig.\ref{fig:C-HiDeNN convergence}(b,d). These two plots show that C-HiDeNN curves are always lower than those of FEM for the same $p$. That is, for the same DOFs, C-HiDeNN is always more accurate than FEM. In other words, to achieve the same level of error, FEM needs more DOFs than C-HiDeNN for the same $p$. When $p=3$, for example, FEM requires 5 to 10 times more DOFs than C-HiDeNN-FEM for a given accuracy (see Fig.\ref{fig:C-HiDeNN convergence}(b,d)).

It is important to note that the convolution operations over the local patch domains are the same as those performed in convolutional neural network (CNN). The larger the patch size $s$, the more neurons connectivity. The dilation parameter $a$ acts as a convolution operator, dictating the feature extraction of the data. Thus, many of the well-developed high-performance computing algorithms utilizing parallel architectures (such as GPU and TPU) in CNN can be integrated with C-HiDeNN. \cite{park2023convolution} discusses in detail on how the parallel programming can accelerate C-HiDeNN computation and demonstrates that it can be as fast as commercial FEM software running on CPU.

\subsubsection{Space Time Finite Element Method with DLDC}
This section proposes a DLDC strategy for solving partial differential equations based on a space-time finite element formulation \cite{hughes1988space,hulbert1990space,wang2017time}. Space-time finite element methods are numerical discretization schemes to predict the spatiotemporal responses of dynamic systems. When used to approximate solutions to partial differential equations, they require assumed values of the underlying governing parameters. For structural problems, these parameters may correspond to material properties such as density, stiffness, and damping capacity. For many applications these assumptions are reasonable, but for many others the distribution of properties throughout a material domain is unknown. For the latter, experimental observation is required for the accurate calibration of model parameters. Consider an elastic one-dimensional bar governed by the following equation.

\begin{equation}
\centering
\label{Eq STFEM 1}
EA\frac{\partial^2u}{\partial x^2}-\rho A\frac{\partial^2u}{\partial t^2}+f=0.
\end{equation}

Here, $t$ is time, $x$ is location, $u$ is displacement, $E$ is elastic stiffness, $A$ is cross-sectional area, $\rho$ is density of the material, and $f$ is external excitation. The discretization of the bar must be done through the space and time dimensions simultaneously, and here it is done with linear quadrilateral elements. Details of the Galerkin formulation are included in the Appendix. The nodes are numbered counting upwards through space in the same order for each time instant, beginning with the first time instant. This results in a highly structured space-time stiffness matrix that will be exploited in the developed DLDC method. For this problem with one spatial dimension, $N$ nodes through space and $T$ nodes through time will result in a total of $N\times T$ degrees of freedom. Besides the enforcement of initial and boundary conditions, the stiffness matrix retains a diagonal structure such that the $N$ equations corresponding to a given time instant only have $3N$ non-zero coefficients. This corresponds to the displacement field at time $t$ being only directly related to the displacement fields at times $t-1$ and $t+1$. For a system discretized with linear quadrilateral elements, a piece of the matrix equation's larger diagonal structure is depicted below with subscripts denoting sub-matrix size and superscripts denoting time step.

\begin{equation}
\centering
\label{Eq STFEM 2}
  \begin{bmatrix}
   \mathbf{A}_{N\times N}&\mathbf{B}_{N\times N}&\mathbf{C}_{N\times N}&0&0\\
   0&\mathbf{A}_{N\times N}&\mathbf{B}_{N\times N}&\mathbf{C}_{N\times N}&0\\
   0&0&\mathbf{A}_{N\times N}&\mathbf{B}_{N\times N}&\mathbf{C}_{N\times N}
   \end{bmatrix}
   \begin{bmatrix}
   \mathbf{u}_{N\times 1}^{t-2}\\
   \mathbf{u}_{N\times 1}^{t-1}\\
   \mathbf{u}_{N\times 1}^{t}\\
   \mathbf{u}_{N\times 1}^{t+1}\\
   \mathbf{u}_{N\times 1}^{t+2}\\
   \end{bmatrix}
   =
   \begin{bmatrix}
   \mathbf{f}_{N\times 1}^{t-1}\\
   \mathbf{f}_{N\times 1}^{t}\\
   \mathbf{f}_{N\times 1}^{t+1}\\
   \end{bmatrix}
\end{equation}

Here, $\bf{A}$, $\bf{B}$, and $\bf{C}$ represent repeated sub-matrices within the global stiffness matrix. For the elastic bar governed by an equation with no damping term, it is noted that $\bf{A}$ and $\bf{C}$ are identical. The solution of this matrix equation can be done sequentially, solving only $N$ simultaneous equations at a time before progressing to the next time instant. These computations can be formulated as an autoregressive neural network. Autoregressive neural networks are sequential feed-forward models used to predict sequence values based on prior observations of the sequence \cite{taskaya2005comparative,triebe2019ar-net}. Generally, this can be represented by the following equation with $ARNN$ denoting the autoregressive neural network as a function of prior observation data, where $\mathbf{X}^t$ and $\bm{\varepsilon}^t$ represent the observation and prediction error respectively at step $t$ in a sequence of data.

\begin{equation}
\centering
\label{Eq STFEM 3}
\mathbf{X}^{t}=ARNN(\mathbf{X}^{t-1},\mathbf{X}^{t-2},\mathbf{X}^{t-3},...,\mathbf{X}^{t-p})+\bm{\varepsilon}^t
\end{equation}

Here, the autoregressive model would be considered of order $p$, since the previous $p$ observations are used to calculate a prediction. Predictions spanning multiple time steps may be achieved by using predicted sequence values as inputs for the prediction of further steps in the sequence. By inspection of the space-time finite element matrix equation \ref{Eq STFEM 2}, the following equivalent second order autoregressive neural network can be constructed.

\begin{equation}
\centering
\label{Eq STFEM 4}
\mathbf{u}^{t+1}
=
ARNN(\mathbf{f}^{t},\mathbf{u}^{t},\mathbf{u}^{t-1})
=
\mathbf{C}^{-1}
\left[
\mathbf{f}^{t} - 
\begin{bmatrix}
\mathbf{A}&\mathbf{B}
\end{bmatrix}
\begin{bmatrix}
\mathbf{u}^{t-1}\\  \mathbf{u}^{t} 
\end{bmatrix}
\right]
\end{equation}

At first glance, this re-framing of the space-time finite element formulation as a neural network provides no apparent advantage. The computations required to solve the partial differential equation remain unchanged. However, recalling that the space-time stiffness matrix elements depend on the coefficients of the partial differential equation, it follows that the neural network provides a functional mapping of these coefficients to predicted sequence values. Sub-matrices, $\bf{A}$, $\bf{B}$, and $\bf{C}$, of the space-time stiffness matrix are functions of the differential equation coefficients according to the finite element assembly. In the context of the elastic bar example, given force and prior displacement fields, the constructed neural network maps the elastic stiffness, cross-sectional area, and density of the bar to the displacement field at the next time step.

\begin{equation}
\centering
\label{Eq STFEM 5}
\mathbf{u}^{t+1}
=
ARNN(E,A,\rho : \mathbf{f}^{t},\mathbf{u}^{t},\mathbf{u}^{t-1})+\bm{\varepsilon}^{t+1}
=
\mathbf{C}^{-1}(E,A,\rho)
\left[
\mathbf{f}^{t} - 
\begin{bmatrix}
\mathbf{A}(E,A,\rho) & \mathbf{B}(E,A,\rho)
\end{bmatrix}
\begin{bmatrix}
\mathbf{u}^{t-1}\\  \mathbf{u}^{t} 
\end{bmatrix}
\right]+\bm{\varepsilon}^{t+1}
\end{equation}

This empowers the automatic learning of the differential equation coefficients from observed spatiotemporal data. If $E$, $A$, and $\rho$ are taken as the trainable parameters of the autoregressive neural network, they can be optimized such that the error term, $\bm{\varepsilon}$, is minimized over all time steps. By embedding the space-time finite element method in the architecture of the neural network, variable time step and spatial mesh sizes can be readily accommodated unlike with black-box neural network implementation. Furthermore, while traditional dense or deep neural network models are prone to over-fitting, much less data is needed to accurately predict system responses with this FEM-informed auto-regressive neural network since the number of trainable parameters is limited to the number of equation coefficients. Additionally, unlike traditional deep learning implementations, extrapolation beyond the training data set to new initial and boundary conditions can be expected since the computations performed during a forward pass retain the form of the underlying numerical method.

To numerically demonstrate the proposed finite element-informed auto-regressive neural network, we consider a spring-mass-damper system governed by the following differential equation.

\begin{equation}
\label{Eq STFEM 6}
m\frac{\partial^2u}{{\partial t}^2}+c\frac{\partial u}{\partial t}+ku=f_0sin(\omega t)
\end{equation}

\begin{figure}[h]
\centering
\begin{subfigure}{.5\textwidth}
  \centering
  \includegraphics[width=.8\linewidth]{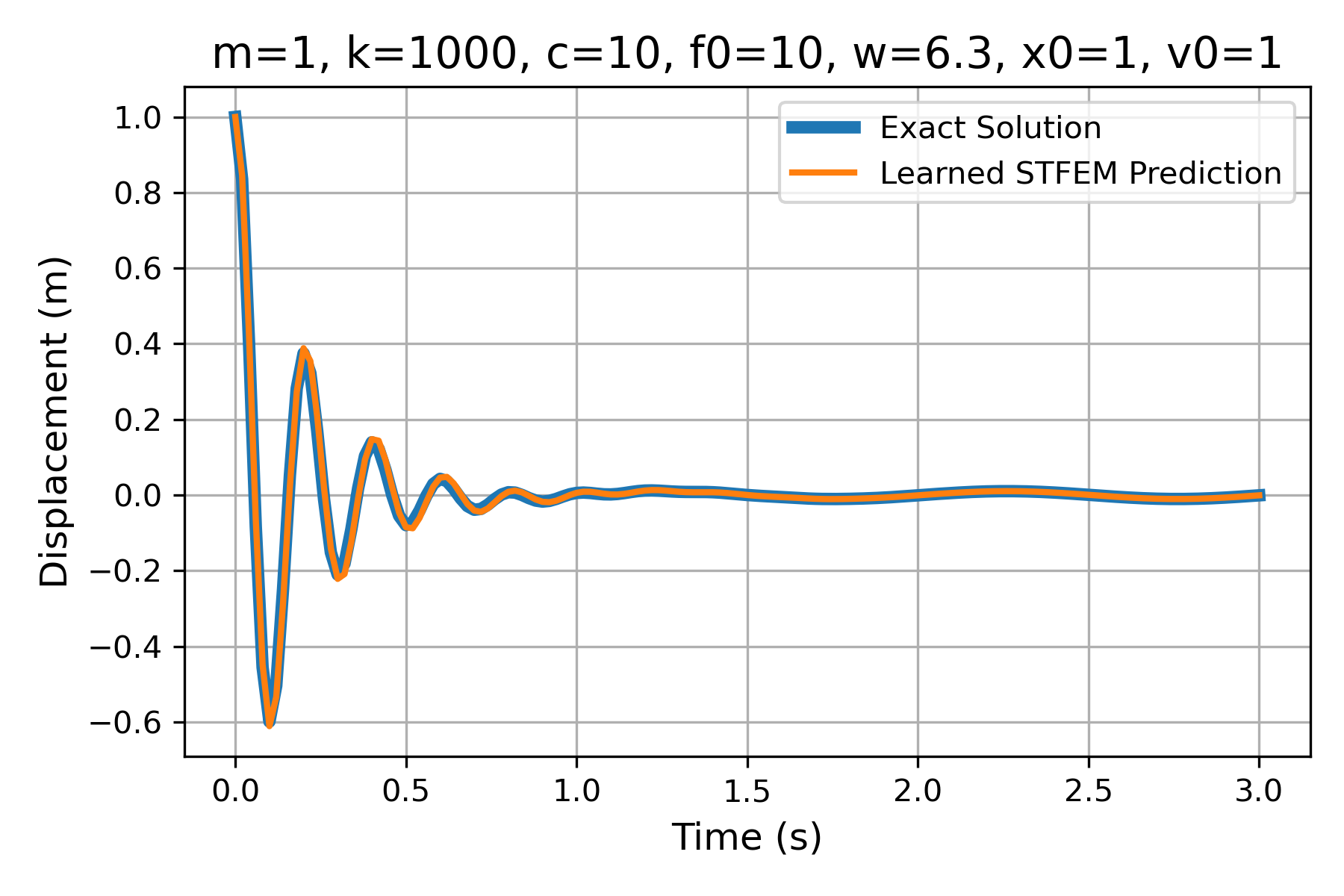}
  \caption{Exact Training Data and Predicted Solution}
  \label{fig:sub1}
\end{subfigure}%
\begin{subfigure}{.5\textwidth}
  \centering
  \includegraphics[width=.8\linewidth]{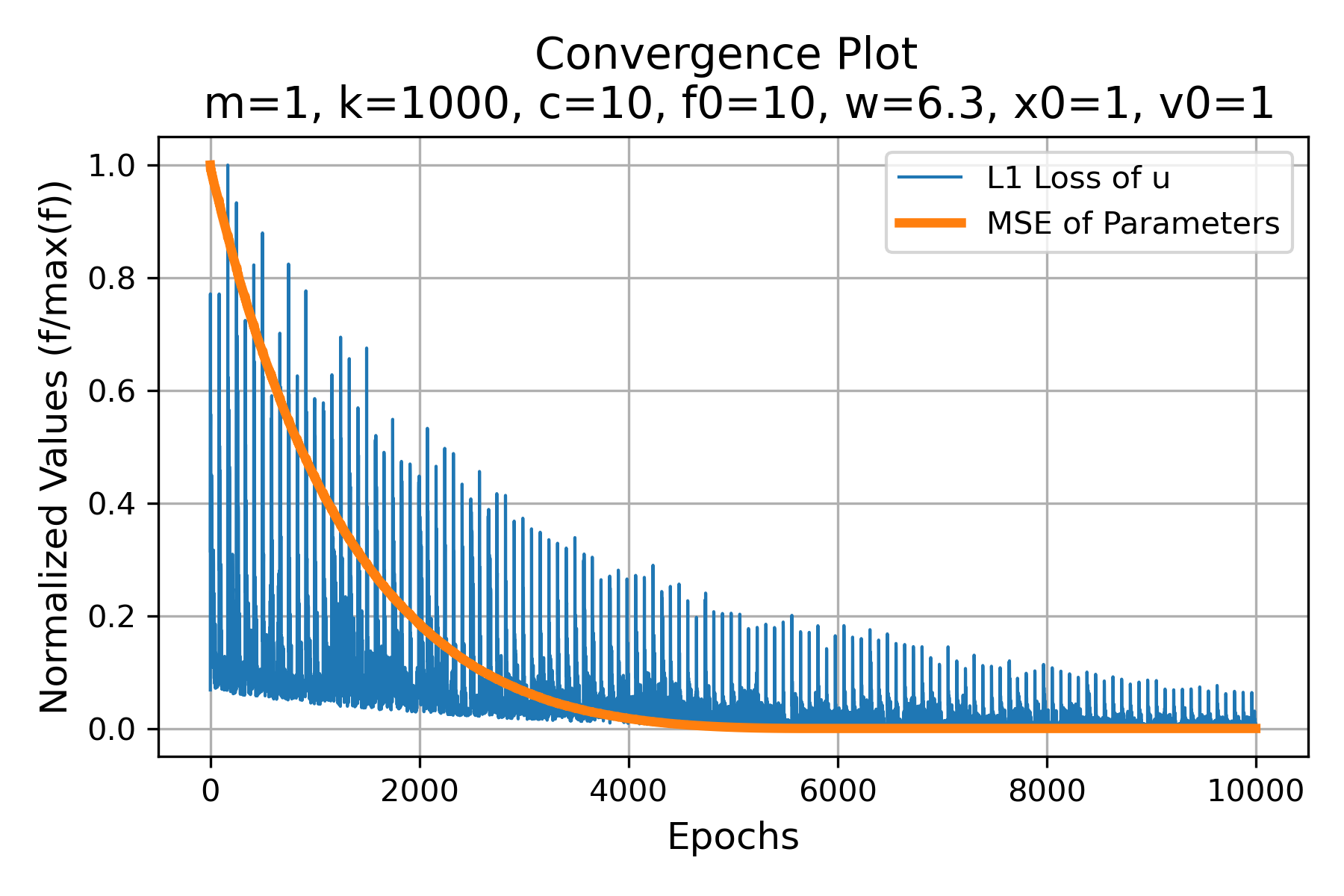}
  \caption{Convergence of the Training Procedure}
  \label{fig:sub2}
\end{subfigure}%
\caption{Training the STFEM Auto-Regressive Neural Network with Exact Data.}
\label{fig:Fig10}
\end{figure}

\begin{figure}[h]
\centering
\begin{subfigure}{.3\textwidth}
  \centering
  \includegraphics[width=1\linewidth]{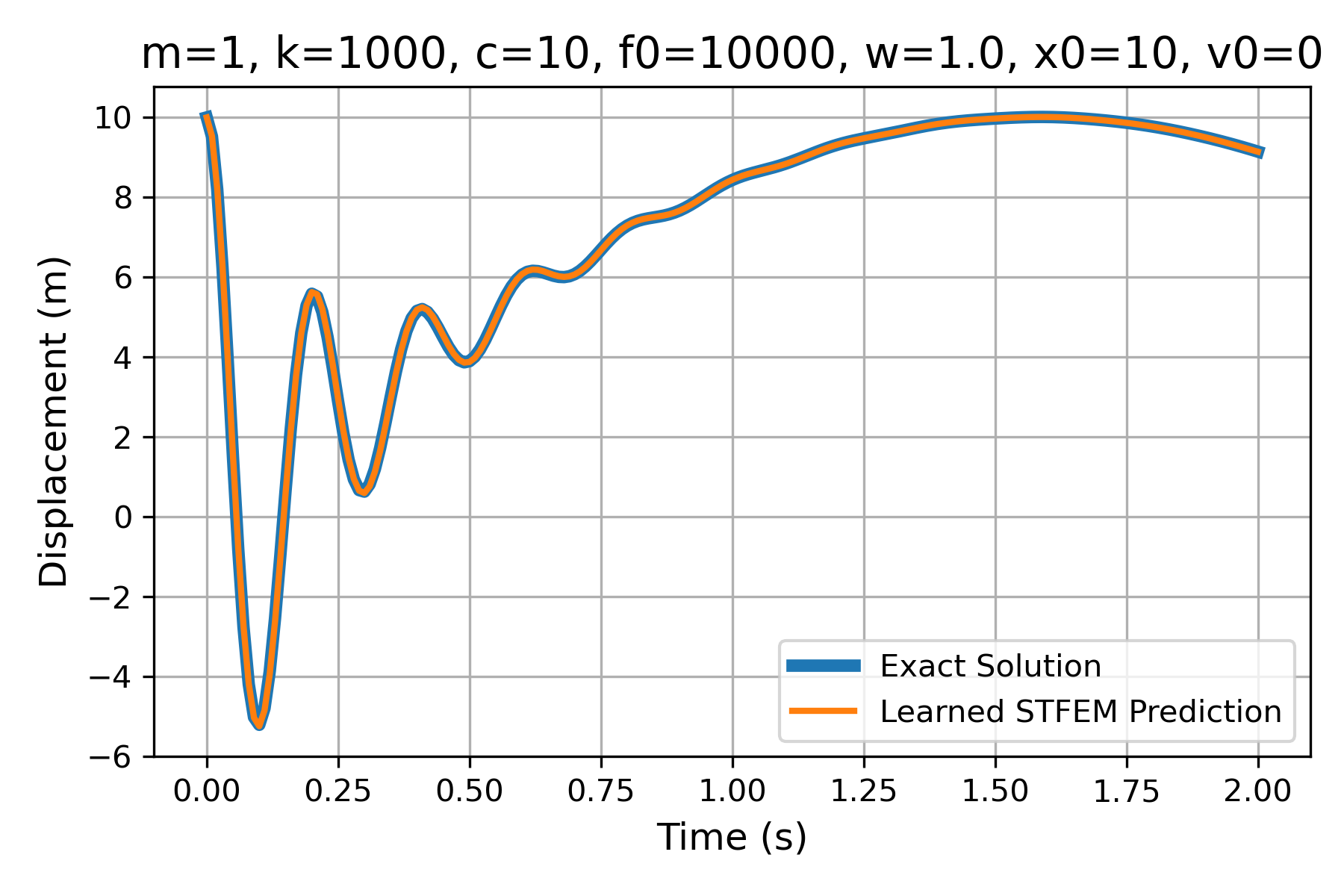}
  \label{fig:sub1}
\end{subfigure}%
\begin{subfigure}{.3\textwidth}
  \centering
  \includegraphics[width=1\linewidth]{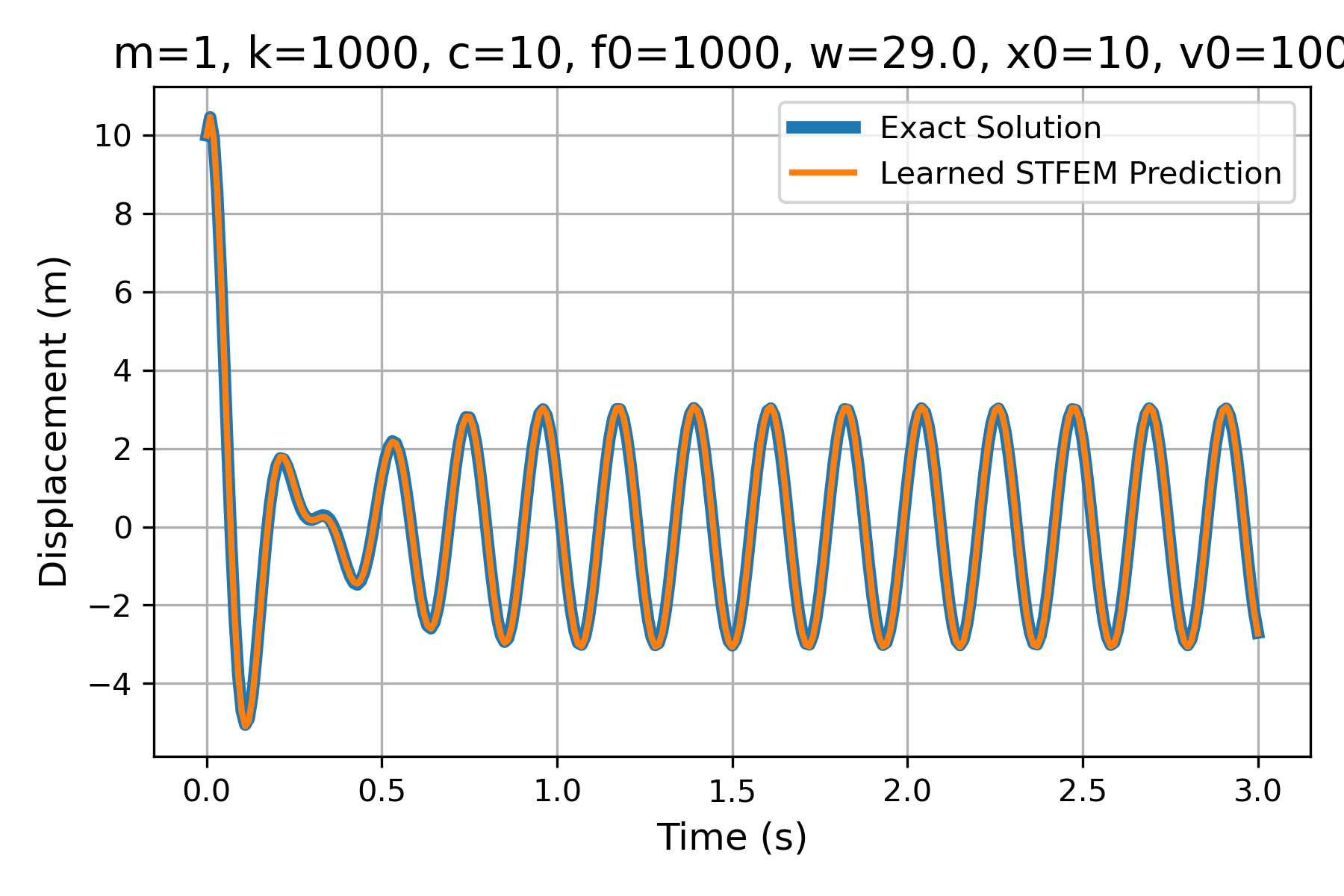}
  \label{fig:sub2}
\end{subfigure}%
\begin{subfigure}{.3\textwidth}
  \centering
  \includegraphics[width=1\linewidth]{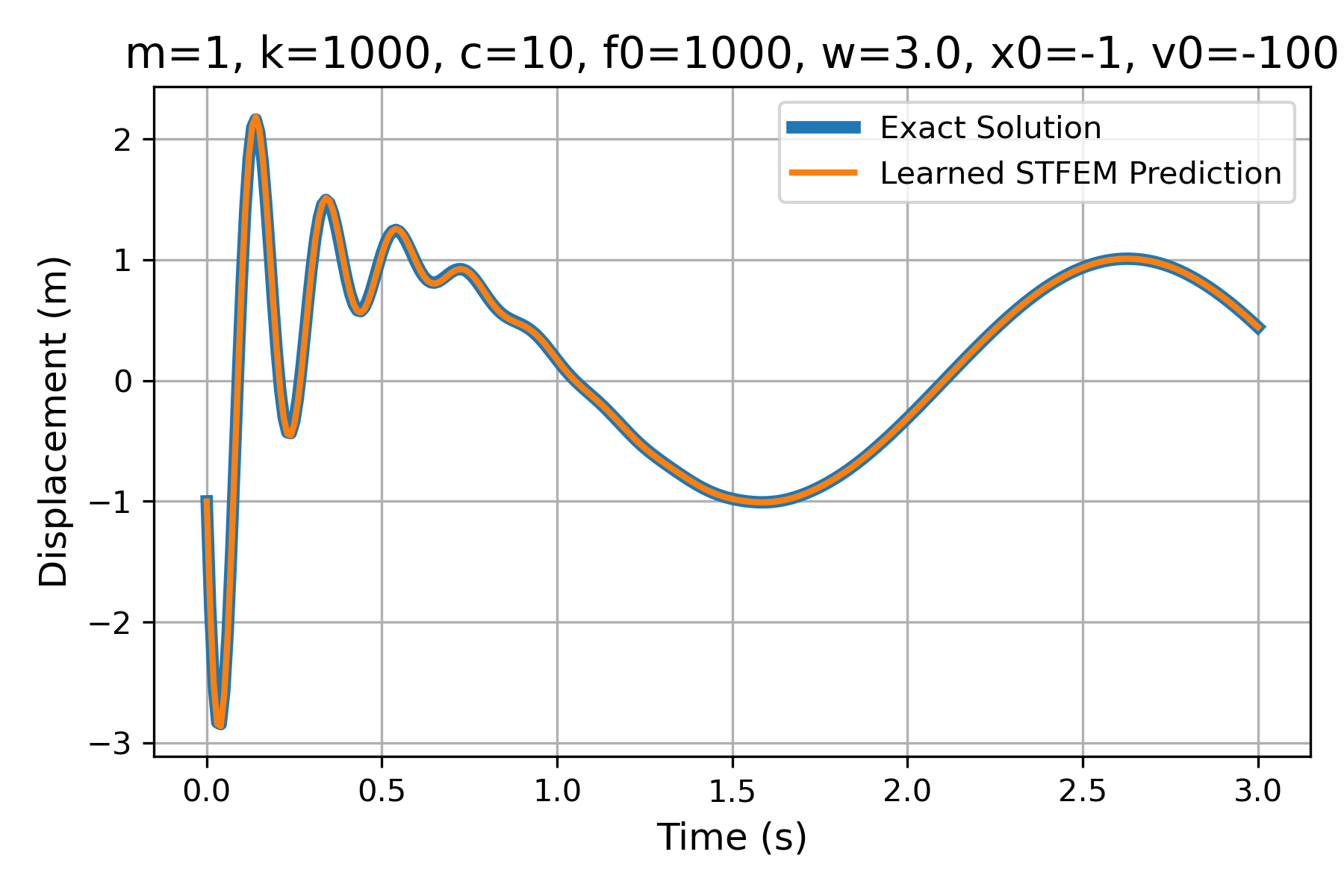}
  \label{fig:sub3}
\end{subfigure}%
\caption{Prediction for New Initial Conditions and Boundary Conditions from Exact Training Data.}
\label{fig:Fig11}
\end{figure}

Here, $t$ is time, $u$ is displacement, $m$ is mass, $c$ is damping coefficient, $k$ is spring stiffness, and $f_0$ and $\omega$ are the maximum magnitude and frequency of an external sinusoidal force. Training data is generated by solving the equation using 150 finite elements spanning 3 time units for the following parameter values and initial conditions: $m = 1, c = 10, k = 100, f_0 = 10, \omega = 2\pi, u(t=0) = 1$, and $\frac{\partial u}{\partial t}(t=0) = 1$. A STFEM autoregressive neural network is implemented with $m$, $c$, and $k$ as the trainable parameters. The sum of the absolute differences, also known as the L1 loss function, between the true displacements, $u^t$, and their predicted values are taken as the objective function to be minimized. The Adam optimization algorithm is used to iteratively update the learnable parameters. A convergence plot of the training procedure is depicted in Fig. \ref{fig:Fig10}(b), which shows that the learnable parameters in the neural network approach the true values of the equation coefficients as the sum of the absolute errors, $|{\varepsilon}^t|$, in the prediction of ${u}^t$ is minimized, resulting in accurate prediction of the equation solution visualized in Fig. \ref{fig:Fig10}(a). The learned values of these coefficients within the neural network can then be used to make predictions for problems with new initial conditions and boundary conditions, generalizing and extrapolating beyond training data like a traditional numerical method while remaining a data-driven machine learning framework. As shown in Fig. \ref{fig:Fig11}, the model is able to accurately predict the behavior of the spring-mass-damper system for a variety of problems with different boundary and initial conditions even though it was only trained on a single case.

\begin{figure}[h]
\centering
\includegraphics[width=0.50\textwidth]{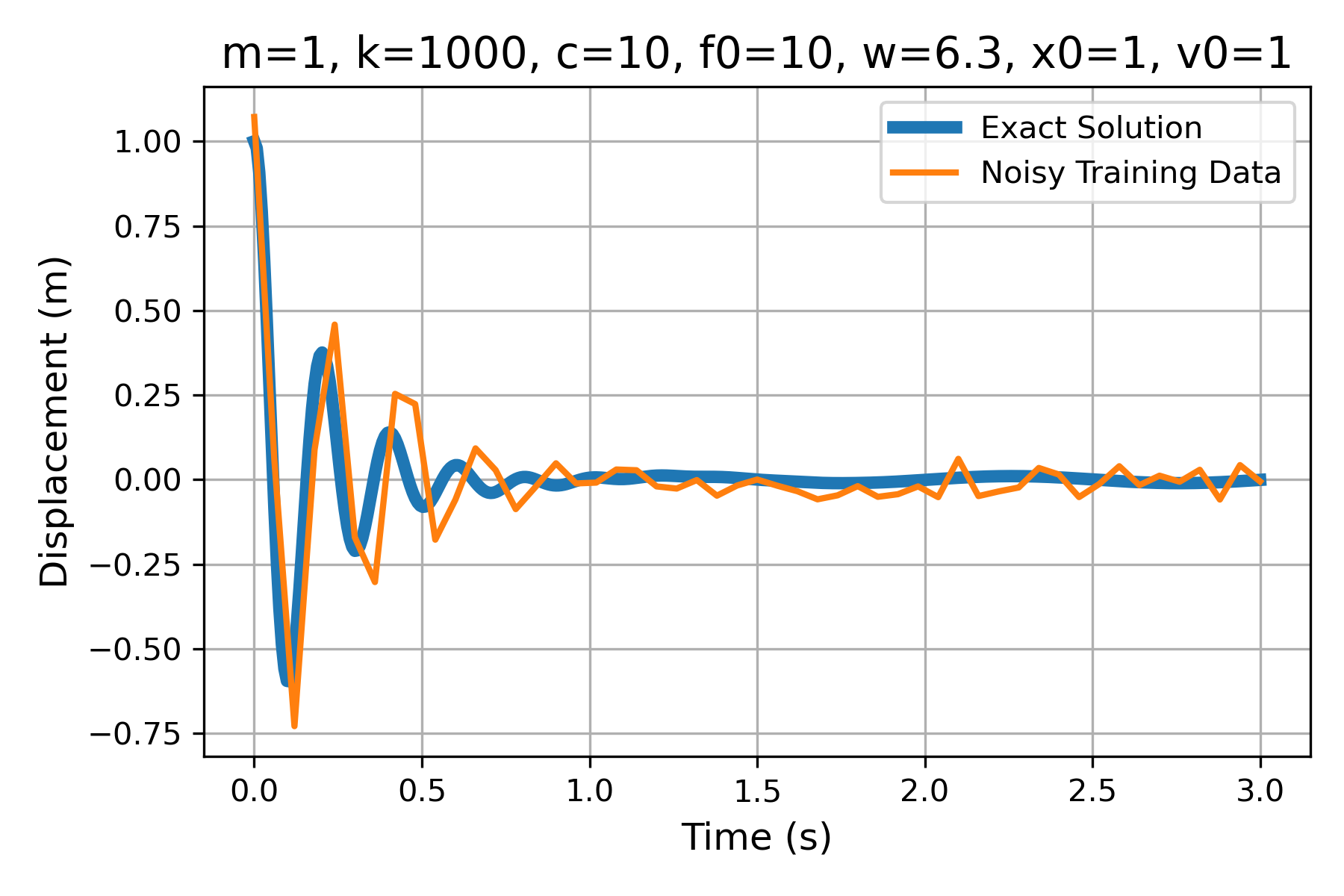}
\caption{True Solution and Generated Noisy Training Data.}
\label{fig:Fig12}
\end{figure}

\begin{figure}[h]
\centering
\begin{subfigure}{.3\textwidth}
  \centering
  \includegraphics[width=1\linewidth]{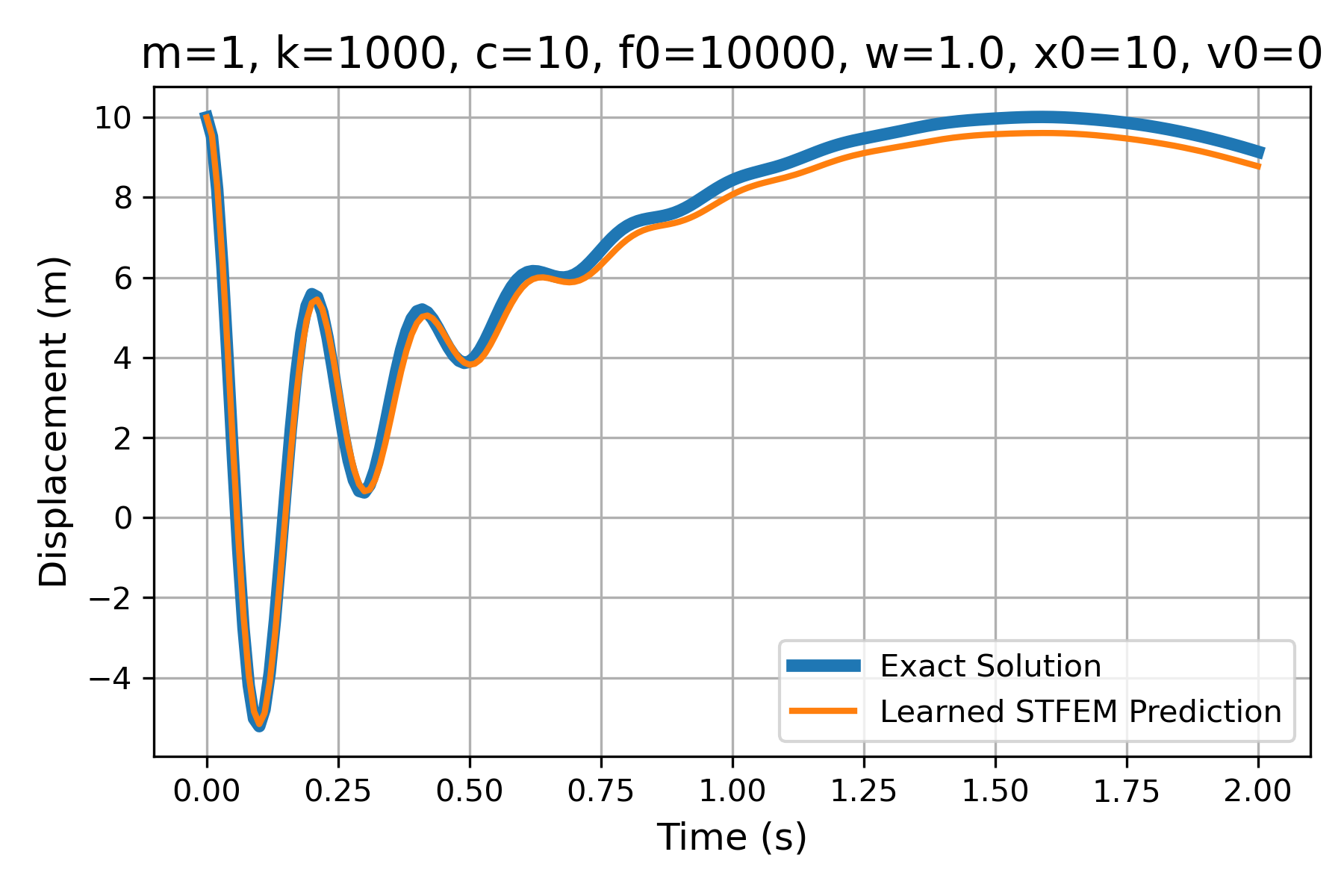}
  \label{fig:sub1}
\end{subfigure}%
\begin{subfigure}{.3\textwidth}
  \centering
  \includegraphics[width=1\linewidth]{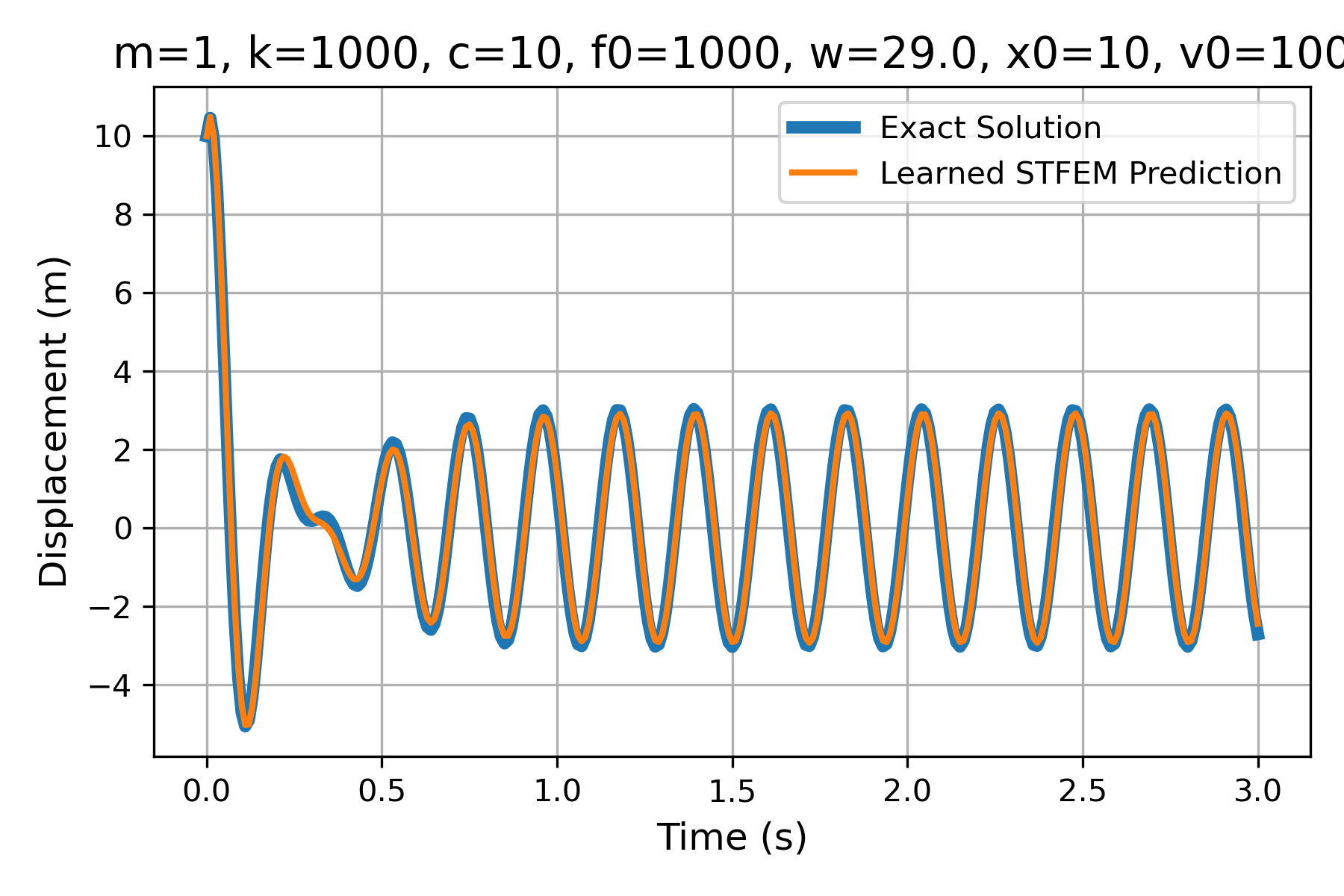}
  \label{fig:sub2}
\end{subfigure}%
\begin{subfigure}{.3\textwidth}
  \centering
  \includegraphics[width=1\linewidth]{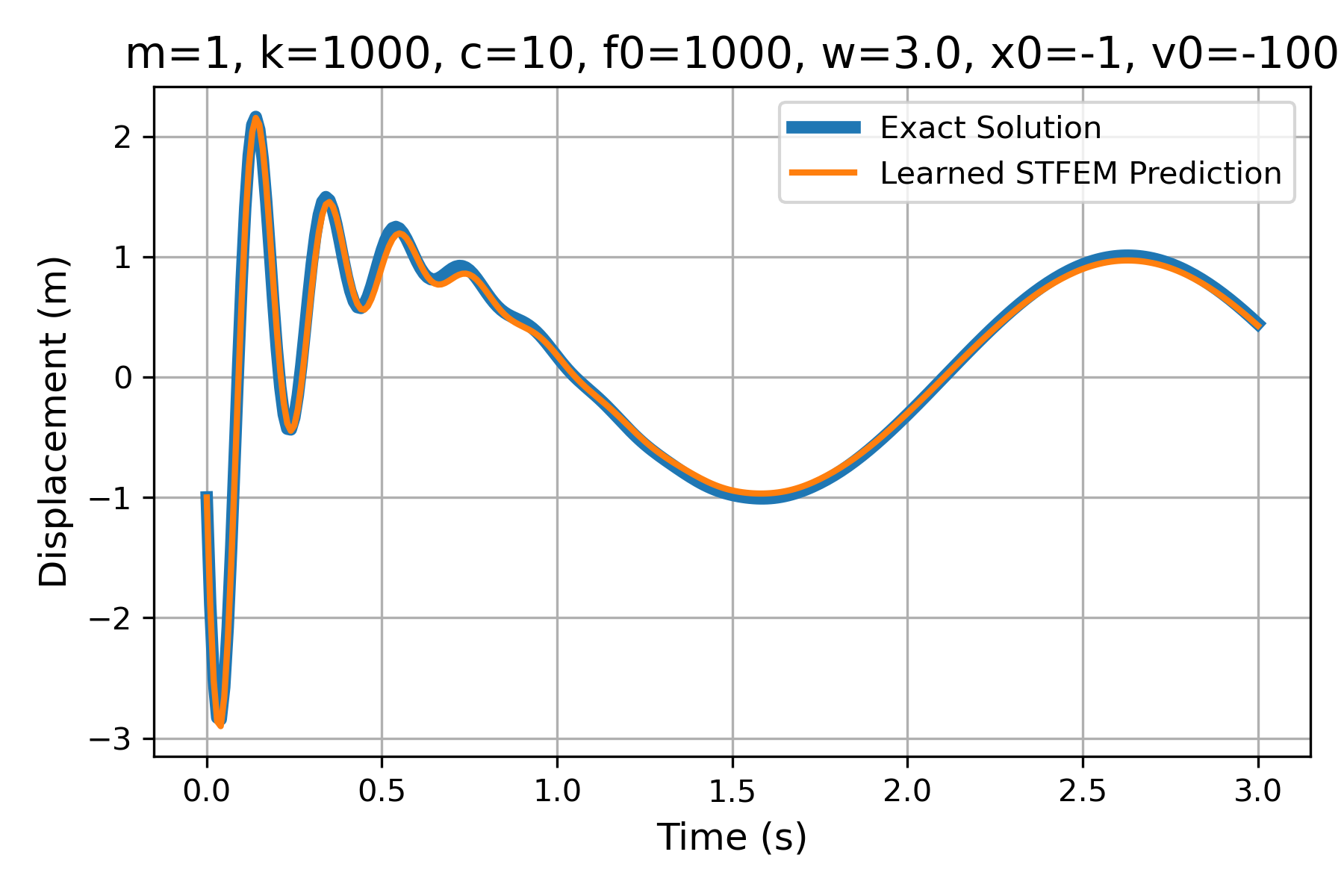}
  \label{fig:sub3}
\end{subfigure}%
\caption{Prediction for New Initial Conditions and Boundary Conditions from Noisy Training Data.}
\label{fig:Fig13}
\end{figure}

Gaussian noise of mean 0 and variance 0.001 is added to a new set of training data for the same problem, which is now solved with a coarser 50 element mesh spanning 3 time units. The resulting time series is shown in Fig. \ref{fig:Fig12}, and Fig. \ref{fig:Fig13} illustrates that even when provided with a very small and noisy training data set the STFEM structured auto-regressive neural network can accurately capture the system behavior across a variety of new conditions. The ability of the STFEM neural network to produce accurate predictions for cases beyond the training data set is not typical of traditional neural networks, which typically only perform well when applied to cases similar to those seen during training.

Dense and deep neural networks offer great function representation capacity, which make them powerful tools for predicting complex phenomena if they are provided with sufficiently large and diverse datasets to learn from. With limited data, neural networks exhibit overfitting due to their large number of trainable parameters leading to poor interpolation capability. Additionally, the extrapolation capability of neural networks to make good predictions beyond the range of their training data is generally very poor with “black box” implementation. These qualities preclude neural networks from solving engineering problems where the collection of data is expensive or cannot span the space of desired prediction capability.
By designing an autoregressive neural network architecture consistent with the space-time finite element method, the generalizability and interpretability of the numerical method can be preserved while harnessing the advantage of data-driven machine learning to reduce the need to make potentially incorrect parameter assumptions. The STFEM neural network may be trained on video data for only a single loading case to produce reasonably accurate parameter values. This intelligently designed neural network, which captures mesh connectivity and spatiotemporal discretization, both reduces the cost of model training and increases robustness to data noise.
The simple neural network architecture, with limited learnable parameters, eliminates the need for collecting a massive quantity of data. In fact, with exact data, only a few time observations are required to eventually converge to the true system parameters, which could then be used to make predictions across a wide variety of boundary conditions in the same manner as finite element analysis. Since the functional structure of the neural network is specified as the finite element method, the collection of additional data only serves to reduce inaccuracies arising due to noise in the observations. This contrasts traditional deep-learning models which seek to learn a functional structure from scratch, hence their massive data requirement.

\subsection{Solving Integral Equations }
The section will discuss how the DLDC methods can be used to solve the integral equations. For brevity and the underlying context of the article, this section will focus on how to solve the Lippmann-Schwinger equation \cite{lippmann1950variational} efficiently with the understanding of the numerical methods through DLDC. This form of integral equation has myriads of applications including scattering \cite{gopal2022accelerated} and multiscale mechanics of materials \cite{zecevic2022new,moulinec1998numerical}. A general integral equation appears as,

\begin{equation} \label{Eq 13}
    f(x) = \int_{a}^{b} \mathcal{K}(x,y)u(y)dy, 
\end{equation}

 where $u$ is an unknown function, $f$ is a known function, and $\mathcal{K}$ is the kernel function. Solution of the integral equation depends on determining the kernel function. The importance of this kernel function will be clear if discussed in the context of solving a differential equation through Green's function. For example, consider the following equation,    

\begin{equation}
    \mathcal{L}u(x) = f(x), \forall x \in \mathcal{\mathbf{R}}
\end{equation}

where, $\mathcal{L}$ is a linear differential operator, $u(x)$ is the unknown function, and $f(x)$ is the forcing function. The solution of this differential equation is given by, 

\begin{equation} \label{Eq 15}
    u(x) = \int \mathcal{G}(x-y)f(y)dy,
\end{equation}

Here, $\mathcal{G}(\cdot)$ is the Green's function, and the following relationship holds true, 

\begin{equation} \label{Eq 16}
    \mathcal{L}\mathcal{G}(x,y) = \delta(x-y),
\end{equation}

where  $\delta(\cdot)$ is the Delta function. If one compares Eq. \ref{Eq 13} and Eq. \ref{Eq 15}, it will be apparent that for linear differential operator $\mathcal{L}(\cdot)$, the Green's function is the kernel function. If one solves for this Green's function, the differential equation will be solved automatically. It is very hard to analytically solve for this Green's function. Usually, this Green's function is determined numerically. 

In recent operator learning paradigm \cite{kovachki2021neural,li2020fourier}, this kernel function is replaced by a convolutional neural network. This neural network or kernel is trained by adding subsequent convolutional layers. In order to achieve resolution independence, in Fourier Neural Network (FNO) the input space of the function variables is taken to Fourier domain and trained \cite{li2020fourier}. Another approach is using graph kernels \cite{li2020neural} where a subset of the domain points from input and the output functional space are considered to construct a graph and convolutional kernel is employed for training. While these networks are shown to solve complex problems including Navier-Stokes equation \cite{li2020fourier}, when the domain to solve becomes very large for example large-scale simulation of 3D microstructure as shown in \cite{saha2021macroscale,saha2021microscale}, the number of training parameters may explode and a lot of offline training efforts are required. In this section, using the concepts of mechanics and DLDC, the article gives formulations for two methods that can reduce the computational burden when solving large-scale engineering problems through convolutional kernels. 

\begin{figure}[h]
\centering
\includegraphics[width=0.85\textwidth]{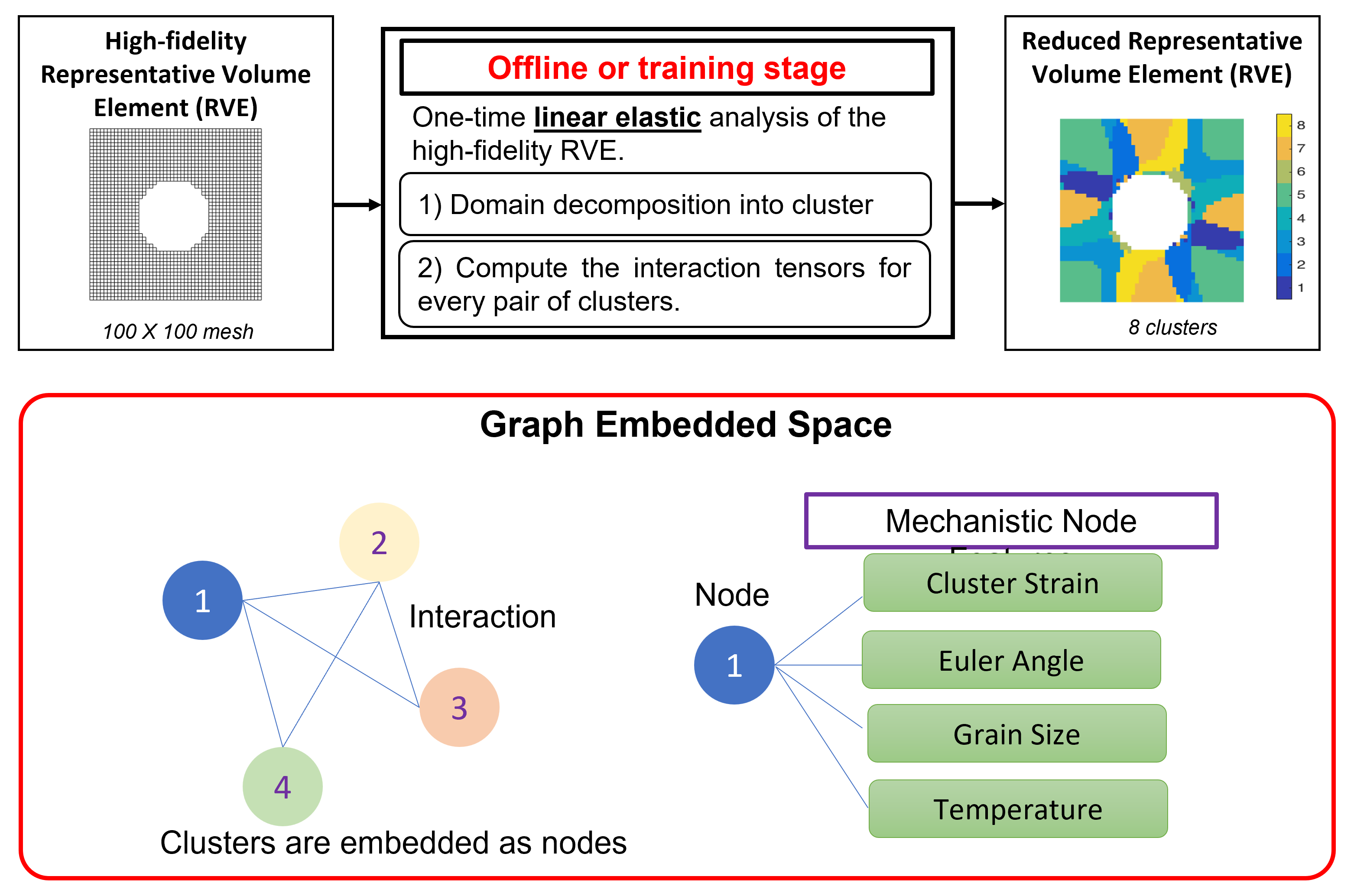}
\caption{The idea of solving a problem with SCA is explained at the top row of the figure. The 2D representative volume element with $100\times100$ mesh is discretized into 8 clusters. In DLDC method, each cluster is assumed to be a node of the graph and mechanistic functional spaces are taken as nodal attributes. The interaction matrix is similar to nodal interaction among the clusters.}
\label{fig:SCA_nodes}
\end{figure}

As a demonstration, the article chooses to solve resort to the self-consistent clustering analysis \cite{liu2016self}. The SCA algorithm has demonstrated considerable success in modeling damage, failure, and mechanical response of multiscale materials system \cite{yu2019self,kafka2021image}. The idea of SCA is to group similar material points together in a structure based on their elastic response, and solve the Lippmann-Schwinger equation only for the group or cluster of material points instead of the entire domain.  More details can be found here \cite{liu2016self}. The following section will touch the basics of SCA for convenience of discussion. 

The homogenization problem for a representative volume element (RVE) can be modeled as the following integral equation. 

\begin{equation} \label{Eq 17}
    \bm{\epsilon(\mathbf{x})} = \bm{\epsilon^{0}(x)}-\int_{\Omega} \mathcal{G}(x,y):\left[\bm{\sigma(y)}-\bm{C^{0}(y):\epsilon(y)}\right]dy,
\end{equation}

where $\bm{\epsilon(x)}$ is local strain at material point $x$, $\Omega$ is the physical domain (RVE), $\bm{\epsilon^{0}(x)}$ is background strain, $\bm{\sigma(y)}$ is the local stress, $C^{0}$ is the reference material stiffness. With SCA, the Eq. \ref{Eq 17}  is not solved for all the material points rather on clustered domain. The clustering is done by K-means clustering \cite{hartigan1979algorithm}. At first, a small, elastic load is applied to the RVE and the resulting strain concentration matrix is stored. Based on this strain concentration matrix, the RVE is clustered (see Figure \ref{fig:SCA_nodes}). On this clustered domain the Eq. \ref{Eq 17} becomes 

\begin{equation} \label{Eq 18}
\Delta\bm{\epsilon^{I}} = \Delta\bm{\epsilon^{0}} - \sum_{j=1}^{k}\left[\frac{1}{c^{I}|\Omega|} \int_{\Omega} \int_{\Omega}\chi^{I}(x)\chi^{J}(y) \mathcal{G}(x,y)dxdy\right]:\left[\Delta\sigma^{I}(y) - C^{0}(y):\Delta\epsilon^{J}(y)\right].   
\end{equation}

In the equation, $c$ is cluster identifier, $k$ is the number of clusters, $I, J$ are the cluster indices, and $\chi^I(x)$ is the cluster variables which assumes value 1 when $x$ is in cluster $I$ and zero otherwise.  If one compares Eqn. \ref{Eq 17} and \ref{Eq 18}, it will be apparent that the Green's function is replaced by the convolution operation also known as the interaction tensor. 

\begin{equation} \label{Eq 19}
    \bm{D}^{IJ} = \frac{1}{c^{I}|\Omega|} \int_{\Omega} \int_{\Omega}\chi^{I}(x)\chi^{J}(y) \mathcal{G}(x,y)dxdy.
\end{equation}
 
In conventional SCA method, this interaction tensor is pre-computed as so-called "offline" database. Based on this pre-computed kernel/interaction tensor, the strain on each cluster is computed from Eq. \ref{Eq 19}. With the DLDC method, the article will show how to determine this kernel on the cluster centroids so that the number of parameters to be trained can be reduced.  

Using graph kernel networks, the Eq. \ref{Eq 18} can be solved on reduced functional space. The justification comes from the mathematical form of the graph kernel neural networks. The graph kernel network has the following form for variable update:

\begin{equation} \label{Eq27}
    v_{t+1}(x) = \sigma\left(Wv_{t}(x) + \int_{\Omega} \mathcal{K}(x,y,a(x),a(y))v_{t}(y)dy\right). 
\end{equation}

Here, $v_{t+1}$ is the output variable in a transformed functional space at $(t+1)$-th time step, $\sigma$ is the activation function, $W$ is weight, $(x,y,a(x),a(y)$ is the edge index in the input function space where $x$ and $y$ are coordinates, and $a(x)$ and $a(y)$ are values of input function. More details can be found here []. The discretized form of the Eq. \ref{Eq27} is:

\begin{equation} \label{Eq28}
    v_{t+1}(x) = \sigma\left(Wv_{t}(x) + \frac{1}{|N(x)|}\sum_{y \in N} \mathcal{K}(x,y,a(x),a(y))v_{t}(y)dy\right). 
\end{equation}

\begin{figure}[h]
\centering
\includegraphics[width=0.95\textwidth]{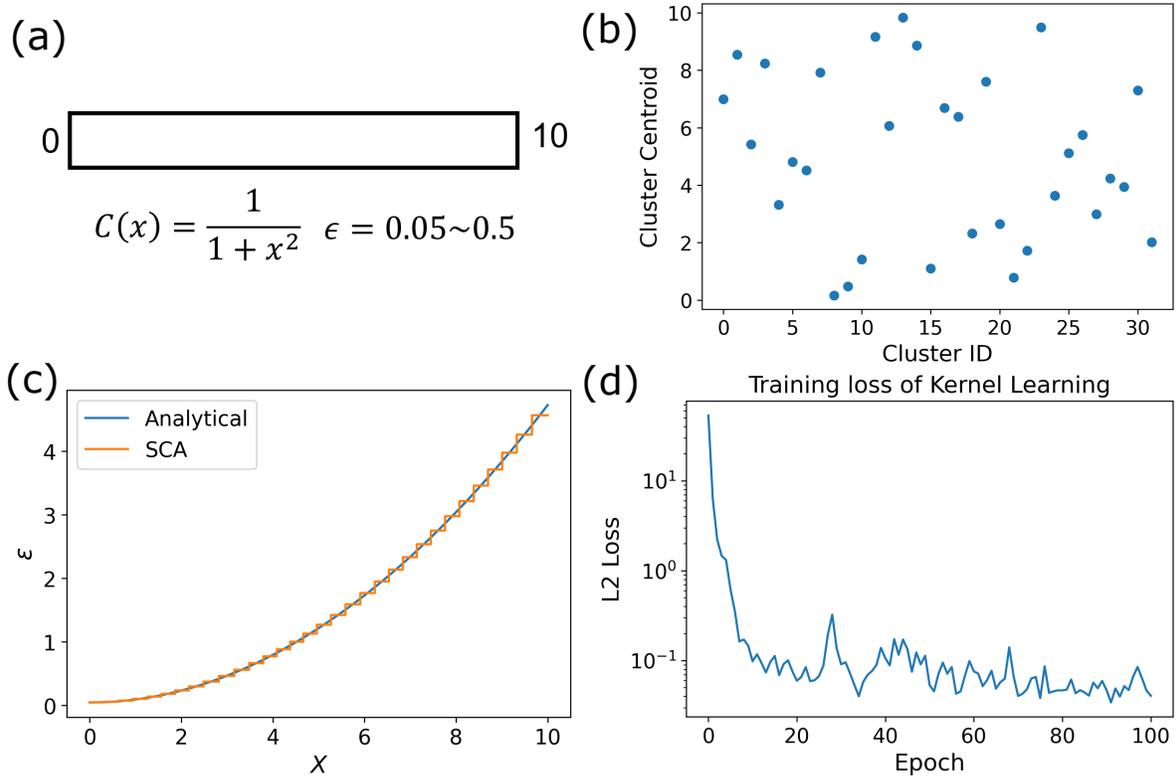}
\caption{Problem statement for the one-dimensional self-consistent clustering analysis. a) The one-dimensional problem domain for training the network, b) The cluster centroids after k-means clustering, c) Comparison between analytical and SCA solution for 33 clusters, d) Evolution of error during training.}
\label{fig:SCA_results}
\end{figure}

\begin{figure}
\centering
\includegraphics[width=0.55\textwidth]{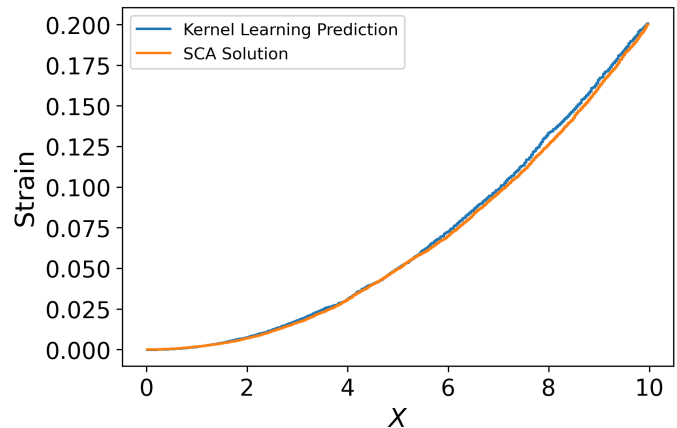}
\caption{Comparison between kernel learning prediction and SCA for 300 clusters.}
\label{fig:SCA_300}
\end{figure}

In this equation, $N(x)$ defines a neighborhood inside radius $r$. Drawing similarities with Eq. \ref{Eq 18}, it is apparent that the interaction tensor and the summation term in Eq. \ref{Eq28} have similarities. The only difference is that, in graph kernel network the neighborhood nodes can be selected by changing the radius of whereas the SCA is inherently global when computing the interaction tensor. As an example, let us consider a one-dimensional domain of length 10 as shown in Figure \ref{fig:SCA_results} (a) discretized into 1000 points in equal intervals. The material is a linear elastic one with stiffness varying as $C(x)=\frac{1}{1+x^2}$. On this domain, the overall applied strain $\epsilon$ is varied from 0.05 to 0.5 and the domain is clustered starting from 2 to 128 clusters. For each such discretization, the cluster centroids are taken as the input samples. Essentially, 1000 material points are reduced 2-128 material points with the assumption that the local strain has a constant value inside each of these clusters. A distribution of possible cluster centroids with 33 clusters is shown in Figure \ref{fig:SCA_results}(b). To generate the data, the solution with SCA is validated against the analytical solution (Figure \ref{fig:SCA_results}(c)). The training is performed with a modified version of the graph kernel network with  6 iterative layers and a radius of 2 (i.e., the interaction is established between only 2 neighbor clusters). The training performance of the network is shown in Figure \ref{fig:SCA_results}. The figure suggests that the training performance is quite good as the normalized mean squared error is going down with number of epochs. A more thorough analysis on the efficacy of such construction is shown in a companion paper \cite{owensahagao2022} where it is shown how varying the hyperparameters of the neural network such as the radius of influence $r$ or number of layers for training. Apart from the reduction of training parameters and domain points, the method can also extrapolate to some extent. It is often identified as the "resolution independence". An example of such extrapolation is shown in Figure \ref{fig:SCA_300}.The figure shows the prediction of the kernel learning algorithm and SCA for 300 clusters and 0.2 applied strain. This particular case is outside the domain of training and testing datasets. In spite of that, we can see that the kernel learning algorithm is as good as SCA. However, further testing is required before claiming that kernel learning method has resolution independence. Nevertheless, the results are promising and shows a way how we can overcome the limitations of the current graph kernel networks while solving three-dimensional engineering problems.

\section{Future Developments}
While the structure of the neural networks proposed in DLDC shows the potential of a general framework of solving fundamental calculus and applied science problems, there are still some limitations that need to be overcome. For example, the forward difference or central difference neural networks presented in the differential calculus section works well only same discretization or sampling intervals are used. This is apparent from the construction of neural networks as well. However, in real-life applications, unstructured data or uneven sampling is quite common. The next step is to extend the DLDC structures similar to differential quadrature method to obtain the derivative of unstructured data. The weights can be learned via either convolutional kernel (similar to C-HiDeNN) or fully-connected neural network or its variant (such as autoencoder). While the DLDC structure gives the exact reproduction of the Gauss quadrature method, it is to be explored if the flexibility of deep neural network can be used to see if it can solve integration that hard to solve with Gauss methods, such as integration around singularities. The authors are currently working to extend C-HiDeNN into a more general form and to solve multiresolution problems such as \cite{mcveigh2006multiresolution,mcveigh2009multiresolution}. While the proposed integral equation solver can significantly reduce the number of training parameters for graph kernels, we still need to train the method with some data. The authors are currently working to extend C-HiDeNN method to solve the integral Lippmann-Schwinger equation so that the convolutional patch function can capture interaction among the cluster centroids.

\section{Conclusions}
The proposed DLDC method aims to combine the rich knowledge of numerical methods with modern data science techniques and create a new perspective on solving challenging problems in engineering and physical sciences. The building blocks of the DLDC are still in the making. The purpose of building the tools is to make a general and easy to follow process of making a deep neural network to solve any governing differential equation. The vision is to combine all the benefits from different numerical methods and propose a unified approach to solve engineering problems. At the same time, the DLDC methods can be used as a teaching tool for calculus in K-12 education system. These tools will simultaneously introduce the students with deep learning and applied calculus. This is expected to increase participation in STEM programs from students. Moreover, the DLDC methods try to avoid multiple training and can be more flexible using the deep neural network. The DLDC version of finite element method i.e., C-HiDeNN shows higher accuracy compared to FEM with same degrees of freedom. The proposed integral equation solver potentially avoids large number of training parameters that was the bottleneck for using graph kernel network. However, further research is required to establish the DLDC as viable method to solve for general problems in science and engineering.

\begin{acknowledgements}

The authors would like to acknowledge the support of National Science Foundation (NSF, USA) grants CMMI-1762035 and CMMI-1934367 and AFOSR, USA grant FA9550-18-1-0381.
C. Park and S. Saha would like to thank the Division of Orthopaedic Surgery and Sports Medicine at Ann and Robert H. Lurie Children's Hospital for their philanthropic grant. The authors would like to acknowledge the contribution of Alberto Ciampaglia, Department of Mechanical and Aerospace Engineering, Politecnico di Torino, Italy and Visiting Researcher, Department of Mechanical Engineering, Northwestern University to this article. 

\end{acknowledgements}

\section{Appendix 1: Space-Time Finite Element Formulation for a 1D Elastic Bar}

The governing equation is:
\begin{equation} 
\label{Eq App 1}
EA\frac{\partial^2u}{\partial x^2}-\rho A\frac{\partial^2u}{\partial t^2}=0.
\end{equation}

For a bar with left end at $x=L$ and right end at $x=R$, vibrating from time $t=0$ to time $t=T$, the above equation can be converted to a weak form by applying the principle of virtual work.
\newcommand{\integral}[4]{\int\limits^{#1}_{#2}\int\limits^{#3}_{#4}}
\begin{equation} 
\label{Eq App 2}
\integral{T}{0}{R}{L} \tilde{v}
\left(EA\frac{\partial^2u}{\partial x^2}-\rho A\frac{\partial^2u}{\partial t^2} \right) 
\partial x \partial t = 0 
\end{equation}

\begin{equation} 
\label{Eq App 3}
\integral{T}{0}{R}{L}
\left(
  EA\frac{\partial \tilde{v}}{\partial x} \frac{\partial^2 u}{\partial x^2}
- EA\frac{\partial \tilde{v}}{\partial x} \frac{\partial u}{\partial x} 
- \rho A\frac{\partial \tilde{v}}{\partial t} \frac{\partial^2 u}{\partial t^2}
+ \rho A\frac{\partial \tilde{v}}{\partial t} \frac{\partial u}{\partial t} 
\right) 
\partial x \partial t  = 0 
\end{equation}

\begin{equation} 
\label{Eq App 4}
\left[  \int_{0}^{T} EA\tilde{v}\frac{\partial u}{\partial x}  \,\partial t  \right]_{x=L}^{x=R}
- \integral{T}{0}{R}{L} EA\frac{\partial \tilde{v}}{\partial x} \frac{\partial u}{\partial x} \partial x \partial t
- \left[  \int_{L}^{R} \rho A\tilde{v}\frac{\partial u}{\partial t}  \,\partial x  \right]_{t=0}^{t=T}
+ \integral{T}{0}{R}{L} \rho A\frac{\partial \tilde{v}}{\partial t} \frac{\partial u}{\partial t} \partial x \partial t  = 0
\end{equation}

Finite element shape function matrices can be used to discretize the simplified weak form. For an element spanning $x=x_0^e$ to $x=x_1^e$ through space and $t=t_0^e$ to $t=t_1^e$ through time, the element matrices correspond to the following expressions.

\begin{equation} 
\label{Eq App 5}
\text{Element Stiffness Matrix  } [K]^e = \integral{t_1^e}{t_0^e}{x_1^e}{x_0^e} EA\frac{\partial \tilde{v}}{\partial x} \frac{\partial u}{\partial x} \partial x \partial t
\end{equation}

\begin{equation} 
\label{Eq App 6}
\text{Element Mass Matrix  } [M]^e = \integral{t_1^e}{t_0^e}{x_1^e}{x_0^e} \rho A\frac{\partial \tilde{v}}{\partial t} \frac{\partial u}{\partial t} \partial x \partial t 
\end{equation}

The global boundary conditions are prescribed via the other terms in Eq. \ref{Eq App 4}.

\begin{equation} 
\label{Eq App 7}
\text{Space Boundary Conditions  } = \left[  \int_{0}^{T} EA\tilde{v}\frac{\partial u}{\partial x}  \,\partial t  \right]_{x=L}
,
\left[  \int_{0}^{T} EA\tilde{v}\frac{\partial u}{\partial x}  \,\partial t  \right]_{x=R}
\end{equation}

\begin{equation} 
\label{Eq App 8}
\text{Time Boundary Conditions  } = \left[  \int_{L}^{R} \rho A\tilde{v}\frac{\partial u}{\partial t}  \,\partial x  \right]_{t=0}
,
\left[  \int_{L}^{R} \rho A\tilde{v}\frac{\partial u}{\partial t}  \,\partial x  \right]_{t=T}
\end{equation}

The time boundary condition at $t=0$ corresponds to the enforcement of initial velocity. The time boundary condition at the final time, $t=T$, can be eliminated from the global matrix equation in favor of equations that enforce initial displacement at $t=0$.

\section{Appendix 2: Time Finite Element Formulation for a Spring Mass Damper System}

The governing equation is:
\begin{equation} 
\label{Eq App2 1}
m\frac{d^2u}{dt^2}+c\frac{du}{dt}+ku = f_0sin(\omega t).
\end{equation}

For a mass oscillating from time $t=0$ to time $t=T$, the above equation can be converted to a weak form by applying the principle of virtual work.
\begin{equation} 
\label{Eq App2 2}
\int_{0}^{T} \tilde{v} \left( m\frac{d^2u}{dt^2} +c\frac{du}{dt} +ku-f_0sin(\omega t) \right)  \,dt = 0
\end{equation}

\begin{equation} 
\label{Eq App2 3}
\left[ m\tilde{v}\frac{du}{dt} \right]_{t=0}^{t=T}
- \int_{0}^{T} m \frac{d\tilde{v}}{dt}\frac{du}{dt} \,dt
+ \int_{0}^{T} c \tilde{v} \frac{du}{dt}\,dt
+ \int_{0}^{T} k \tilde{v} u \,dt
- \int_{0}^{T} \tilde{v} f_0sin(\omega t) \,dt
= 0
\end{equation}

Finite element shape function matrices can be used to discretize the simplified weak form. For an element spanning $t=t_0^e$ to $t=t_1^e$ through time, the element matrices correspond to the following expressions.

\begin{equation} 
\label{Eq App2 4}
\text{Element Mass Matrix  } [M]^e = \int_{t_0^e}^{t_1^e} m \frac{d\tilde{v}}{dt}\frac{du}{dt} \,dt
\end{equation}

\begin{equation} 
\label{Eq App2 5}
\text{Element Damping Matrix  } [C]^e = \int_{t_0^e}^{t_1^e} c \tilde{v} \frac{du}{dt}\,dt
\end{equation}

\begin{equation} 
\label{Eq App2 6}
\text{Element Stiffness Matrix  } [K]^e = \int_{t_0^e}^{t_1^e} k \tilde{v} u \,dt
\end{equation}

\begin{equation} 
\label{Eq App2 7}
\text{Element Force Matrix  } [F]^e = \int_{t_0^e}^{t_1^e} \tilde{v} f_0sin(\omega t) \,dt
\end{equation}

The element force matrices account for the force boundary condition due to the external sinusoidal excitation of the sprung mass. The temporal boundary conditions are prescribed via the other term in Eq. \ref{Eq App2 3}.

\begin{equation} 
\label{Eq App2 6}
\text{Time Boundary Conditions  } = \left[ m\tilde{v}\frac{du}{dt} \right]_0
,
\left[ m\tilde{v}\frac{du}{dt} \right]_T
\end{equation}

The time boundary condition at $t=0$ corresponds to the enforcement of initial velocity. The time boundary condition at the final time, $t=T$, can be eliminated from the global matrix equation in favor of equations that enforce initial displacement at $t=0$.

\bibliography{references.bib}

\end{document}